\documentclass[11pt]{article}

\usepackage[utf8]{inputenc}
\usepackage[T1]{fontenc}
\usepackage{textcomp}
\usepackage{amsfonts}
\usepackage[dvipsnames]{xcolor}
\usepackage{amsmath,amssymb,comment}
\usepackage{tikz}
\usepackage{pgfplots}
\pgfplotsset{compat=1.16}
 \usepackage{hyperref}
 \hypersetup{%
  colorlinks=true,
  linkcolor=black,
  citecolor=black,
  urlcolor=blue,
  pdftitle={generalized solutions for a viscous wave model},
  pdfauthor={Claes, Lankeit, Winkler},
  pdfkeywords={},
  bookmarksopen=false,
}

\newtheorem{theo}{Theorem}[section]
\newtheorem{lem}[theo]{Lemma}

\newtheorem{prop}[theo]{Proposition}

\newtheorem{defi}[theo]{Definition}

\newcommand{\mysection}[1]{\section{#1} \setcounter{equation}{0}}
\newcommand{\proof}{{\sc Proof.} \quad}
\newcommand{\proofc}{{\sc Proof} \ }
\newcommand{\be}{\begin{equation} \label}
\newcommand{\ee}{\end{equation}}
\newcommand{\bea}{\begin{eqnarray}\label}
\newcommand{\eea}{\end{eqnarray}}
\newcommand{\bas}{\begin{eqnarray*}}
\newcommand{\eas}{\end{eqnarray*}}
\newcommand{\bit}{\begin{itemize}}
\newcommand{\eit}{\end{itemize}}
\newcommand{\qed}{\hfill$\Box$ \vskip.2cm}
\newcommand{\nn}{\nonumber}
\newcommand{\R}{\mathbb{R}}
\newcommand{\N}{\mathbb{N}}
\newcommand{\pO}{\partial\Omega}

\newcommand{\eps}{\varepsilon}

\newcommand{\supp}{{\rm supp} \, }

\newcommand{\wto}{\rightharpoonup}

\newcommand{\hra}{\hookrightarrow}
\newcommand{\io}{\int_\Omega}
\newcommand{\na}{\nabla}
\newcommand{\Del}{\Delta}
\newcommand{\del}{\delta}
\newcommand{\al}{\alpha}
\newcommand{\lam}{\lambda}

\newcommand{\sig}{\sigma}
\newcommand{\pa}{\partial}
\newcommand{\bom}{\overline{\Omega}}
\newcommand{\Om}{\Omega}
\newcommand{\om}{\omega}

\newcommand{\mult}{\otimes}
\newcommand{\omult}{\odot}

\newcommand{\wh}{\widehat}

\newcommand{\hs}{\hspace*}

\newcommand{\vp}{\varphi}
\newcommand{\vt}{\vartheta}
\newcommand{\lbal}{\left\{ \begin{array}{l}}
\newcommand{\lball}{\left\{ \begin{array}{ll}}
\newcommand{\ear}{\end{array} \right.}

\newcommand{\jl}[1]{{\color{BurntOrange}#1}}

\newcommand{\abs}{\\[5pt]}

\newcommand{\adb}{\allowdisplaybreaks}
\newcommand{\tm}{T_{max}}

\newcommand{\tme}{T_{max,\eps}}

\newcommand{\lan}{\langle}
\newcommand{\ran}{\rangle}
\renewcommand{\div}{{\rm div} \,}

\newcommand{\ueps}{u_\eps}
\newcommand{\veps}{v_\eps}
\newcommand{\Teps}{\Theta_\eps}

\newcommand{\heps}{h_\eps}

\newcommand{\yeps}{y_\eps}

\newcommand{\F}{{\mathcal{F}}}
\newcommand{\Feps}{{\mathcal{F}_\eps}}
\newcommand{\D}{{\mathcal{D}}}
\newcommand{\Rest}{{\mathcal{R}}}

\newcommand{\nas}{\na^s}
\newcommand{\Beta}{B}

\newcommand{\tops}[2]{\texorpdfstring{#1}{#2}}
\newcommand{\norm}[2][]{\|#2\|_{#1}}
\newcommand{\Lom}[1]{L^{#1}(\Omega)}
\newcommand{\f}[2]{\frac{#1}{#2}}
\newcommand{\Lloc}[1]{L^{#1}_{loc}}

\newcommand{\set}[1]{\{#1\}}
\newcommand{\bdry}{|_{\partial\Omega}}

\usepackage{newunicodechar}
\newunicodechar{ℝ}{\mathbb{R}}
\newunicodechar{ψ}{\psi}
\newunicodechar{Δ}{\Delta}
\newunicodechar{ε}{\eps}
\newunicodechar{ℕ}{\mathbb{N}}
\newunicodechar{ξ}{\xi}
\newunicodechar{∞}{\infty}
\newunicodechar{σ}{\sigma}
\newunicodechar{φ}{\varphi}
\newunicodechar{Θ}{\Theta}
\newunicodechar{χ}{\chi}
\newunicodechar{γ}{\gamma}
\newunicodechar{Γ}{\Gamma}
\newunicodechar{ζ}{\zeta}
\newunicodechar{∇}{\nabla}
\newunicodechar{κ}{\kappa}
\newunicodechar{λ}{\lambda}
\newunicodechar{μ}{\mu}
\newunicodechar{β}{\beta}
\newunicodechar{δ}{\delta}
\newunicodechar{α}{\alpha}
\newunicodechar{τ}{\tau}
\newunicodechar{ρ}{\rho}

\begin{document}
\adb
%
%
\title{A model for heat generation by acoustic waves in piezoelectric materials: Global large-data solutions}
\author{
Leander Claes\footnote{claes@emt.uni-paderborn.de}\\
{\small Universität Paderborn}\\
{\small Institut für Elektrotechnik und Informationstechnik}\\
{\small 33098 Paderborn, Germany}
\and
 Johannes Lankeit\footnote{lankeit@ifam.uni-hannover.de}\\
 {\small Leibniz Universität Hannover (LUH),}\\
 {\small Institut für Angewandte Mathematik}\\
 {\small 30161 Hannover, Germany}\\
 {\small and }\\
 {\small Cluster of Excellence PhoenixD, LUH}
\and
Michael Winkler\footnote{michael.winkler@math.uni-paderborn.de}\\
{\small Universit\"at Paderborn}\\
{\small Institut f\"ur Mathematik}\\
{\small 33098 Paderborn, Germany} }
\date{}
\maketitle
\begin{abstract}
\noindent 
A model for the generation of heat due to
mechanical losses during acoustic wave propagation in a solid is considered in a Kelvin-Voigt type framework.
In contrast to previous studies on related thermoviscoelastic models, in line with recent experimental findings
the present manuscript focuses on situations in which elastic parameters depend on temperature.
Despite an apparent loss of mathematically favorable structural properties thereby encountered,
in the framework of a suitably generalized concept of solvability a result on global existence of solutions
is derived under mild assumptions which, in particular, do not involve any smallness condition on the initial data.\abs
%
\end{abstract}
\newpage
\section{Introduction}\label{intro}
As a mathematical subject, the theory of elasticity goes back to Cauchy and has again flourished since the mid-twentieth century (cf.~eg.~\cite{truesdell_cauchy}). Thermoviscoelasticity with its strong nonlinear coupling between heat generation and the deformation of materials offers mathematical challenges at various levels.
Classical one-dimensional existence results include local \cite{racke} or local and, under smallness conditions on initial data, global existence of classical solutions  \cite{slemrod,dafermos_hsiao_smooth,kim,jiang1990} -- although, on the other hand, for large data, $C^2$ solutions may blow up, \cite{dafermos_hsiao,mcintire} -- or global existence of weak solutions, \cite{racke_zheng,cieslak_galic_muha}. Also the long-time behavior of solutions has been studied, \cite{racke_zheng,qin,bies_cieslak}.
Only substantially later than \cite{slemrod,dafermos_hsiao_smooth} were results for higher-dimensional systems obtained, e.g.  \cite{racke90,shibata,blanchard_guibe} or \cite{gawinecki03} (the latter under assumptions of radial symmetry).
Once again, there are nonlinearities that can cause solutions to lose $C^2$ regularity, \cite{racke_bu}.
Existence results concern local-in-time solutions \cite{jiang_racke90,bonetti_bonfanti} or solvability under smallness conditions \cite{racke90,shibata}.\abs
Global weak solutions in three-dimensional domains were found in \cite{owczarek_wielgos} for a system
with irreversible mechanical deformations manifesting in an inelastic constitutive relation taking the form of an ODE.
Without this inelastic effect, heat capacity growing with temperature has been leveraged in several constructions of certain solutions,
\cite{roubicek,blanchard_guibe,pawlow_zajaczkowski_cpaa17,gawinecki_zajaczkowski_cpaa,gawinecki_zajaczkowski}, 
where this growth was relied on
to deal with analytic difficulties stemming from the effects of temperature dilation.
For a simplified system neglecting viscosity-driven transfer of mechanical energy into heat, a statement on global solvability has
recently been derived in \cite{cieslak_muha_trifunovic}.
Further extended models feature fourth-order terms  \cite{yoshikawa_pawlow_zajaczkowski,rossi_roubicek_interfaces13,yoshikawa_pawlow_zajaczkowski_SIMA,roubicek_nodea13,mielke_roubicek} (which may be beneficial in the derivation of a priori bounds) or  additional 'internal variables' describing e.g. plasticity  \cite{roubicek_SIMA10,bartels_roubicek,bartels_roubicek_m2an11,bartels_roubicek_numpde13,paoli_petrov,paoli_petrov_gamm,roubicek_dcdss13}; for related settings regarding adhesive contact see \cite{rossi_roubicek,rossi_roubicek_interfaces13}. \abs
Common to all these precedent studies on higher-dimensional cases seems a concentration on situations in which crucial system
ingredients such as the elastic parameters and viscosities are constants;
recent experimental observations
concerned with the behavior of certain piezoceramics, however, have revealed some partially significant
dependencies of these constituents on temperature.
The present manuscript investigates a multi-dimensional model for thermoviscoelasticity capable of taking such effects into account,
and intends to develop a basic theory of global solvability on the basis of a generalized solution concept that
seems novel in this context; in particular, the design of this framework will be motivated by the ambition to
cope with an apparent loss of some favorable structural properties going along with such modifications.\abs
Before introducing the mathematical main results, let us consider the model with its physical background.\abs
{\bf Application and modeling background.} \quad
The generation of excess heat has a detrimental effect in many industrial and scientific applications.
Not only are involuntary thermal emissions an indicator for a lack of efficiency, the increased temperature may also damage components.
Among the materials that are particularly susceptible to temperature related damage are piezoelectric ceramics~\cite{Rupitsch2019}.
They are used as electromechanical transducers to generate and detect mechanical vibrations and acoustic waves in a variety of applications, ranging from microphones and loudspeakers to ultrasonic welding.
The piezoelectric effect in these ceramics, which also show ferroelectric properties, is only present if the material is polarised.
This intrinsic polarisation vanishes if a certain temperature threshold, the Curie temperature, is exceeded, rendering the piezoelectric ceramic inert.
This is especially problematic in high-power applications, such as ultrasonic bonding and welding, where piezoelectric ceramics are used as actors and are thus a primary source of heat~\cite{Lesieutre1996,Zheng1996}.
Additionally, the resonant behavior, which is crucial for the function of these devices, shows strong dependence on temperature~\cite{Wellendorf2023}.
Therefore, it is of utmost importance to accurately consider the heat generated by the acoustic waves in and around the piezoelectric material.\abs
To account for mechanical losses		
during acoustic wave propagation, a suitable refinement of Hooke's law is required.
Known as the apparently most basic description for the linear mechanical behavior of a solid,
Hooke's law postulates the mechanical stress tensor \(T = {(T_{ij})}_{i,j \in \{1, 2, 3\}}\in \R^{3 \times 3}\) to be related to the mechanical strain tensor \(S = {(S_{kl})}_{k,l \in \{1, 2, 3\}}\in \R^{3 \times 3}\) through the forth-rank elasticity tensor \(C = {(C_{ijkl})}_{i,j,k,l \in \{1, 2, 3\}}\in \R^{3 \times 3 \times 3 \times 3}\)~\cite{Nye1985} according to the tensor product
\be{stressstrainsymmetric}
	T=C:S,
\ee
that is, to the relation
\bas
	T_{ij} = \sum_{k,l=1}^3 C_{ijkl} S_{kl},\qquad\qquad i,j\in\set{1,2,3} ,
\eas
(cf.~also below for a more compact summary on notation used here);
we note that
both the stress and strain tensor are symmetric in the sense of satisfying $T_{ij} = T_{ji}$ and $S_{ij} = S_{ji}$ for
$i,j\in\{1,2,3\}$.
Now in further development of this, the Kelvin-Voigt model proposes to describe
mechanical losses, and thus the generation of heat, by means of the modified relation (\cite{GutierrezLemini2014})
\begin{equation}
	T = C : (S + \tau S_{t}),
	\label{eq:kelvin-voight}
\end{equation}
where \(\tau \in \R \) is the retardation time constant quantifying losses and \(S_{t}\) denotes the time derivative of the strain~\cite{Meyers2008};
thus containing only one single additional constant compared to the above purely elastic model, this Kelvin-Voigt law is one of the fundamental viscoelastic material models.
The differential equation for the displacement field \(u = {(u_{i})}_{i \in \{1, 2, 3\}} \in \R^{3}\) in a material described by the Kelvin-Voigt model can be derived from the Cauchy momentum equation~\cite{Boley2012}
\begin{align*}
	u_{tt} &= \frac{1}{\rho} \sum_{j=1}^3 \partial_j T_{\cdot j}
	         = \frac{1}{\rho}\div T ,
\end{align*}
and the definition of the strain via the symmetric gradient
\begin{align}
	S_{kl} &= \frac{1}{2}(\partial_l u_k + \partial_k u_l)
				 = (\nas u)_{kl},\qquad\qquad k,l\in\{1,2,3\}.
	\label{eq:strain_def}
\end{align}
With $\rho>0$ denoting the density of the material, the resulting viscous wave equation
\begin{equation}
	\rho  u_{tt} = \div (C : \nas u) + \tau \div (C : \nas u_{t})
	\label{eq:visc_wave}
\end{equation}
does not only model the behavior of a viscoelastic solid described by the Kelvin-Voigt model but also arises in fluid acoustics from a linearized form of the Navier-Stokes equations~\cite{Rudenko1977}.
For \(\tau = 0\) and scalar quantities, (\ref{eq:visc_wave}) takes the form of a classical wave equation with phase velocity \(c_\mathrm{ph} = \sqrt{C / \rho}\).
The mixed third-order derivative term quantifies absorption for acoustic waves and is directly related to the viscosity when considering fluids~\cite{Rudenko1977};
this term thus describes a conversion of mechanical energy into thermal energy.\abs
To quantify the energy conversion process, the work \(P \in \R \) done on a solid can be determined by the scalar product of the mechanical stress \(T\) with the time derivative of the mechanical strain \(S_{t}\)~\cite{Boley2012} according to
\begin{align*}
	P &= \sum_{i,j=1}^3 T_{ij} S_{ij,t} 
		= \langle T, S_{t} \rangle ,
\end{align*}
which in conjunction with \eqref{eq:kelvin-voight} yields the identity
\begin{equation*}
	P = \langle C : S, S_{t} \rangle + \langle \tau C : S_{t}, S_{t} \rangle .
\end{equation*}
The first term in this expression describes reversible energy storage due to elastic deformations;
for harmonic processes especially, it is easily shown that the time average of \(\langle C : S, S_{t} \rangle \) is zero.
The second term quantifies the conversion of mechanical work into thermal energy.
It can be described in terms of the strain \(S\) or, using (\ref{eq:strain_def}), in terms of the displacement \(u\)~\cite{Tauchert1967,Boley2012},
\begin{align*}
	Q = \langle \tau C : S_{t}, S_{t} \rangle
	= \langle \tau C : \nas u_{t}, \nas u_{t} \rangle .
\end{align*}
The quantity \(Q \in \R \) has the physical unit of an energy source density and can be inserted directly as a source for heat generation~\cite{Boley2012} in the parabolic equation
\begin{align}
	c \rho \Theta_t = \lambda \Delta \Theta + \langle \tau C : \nas u_{t}, \nas u_{t} \rangle ,
	\label{eq:thermo_coupled}
\end{align}
for the temperature distribution $\Theta$, where $c>0$ and $\lambda>0$ denote the heat capacity and the thermal conductivity
of the material, respectively; in this sense, (\ref{eq:thermo_coupled}) thus couples thermal and mechanical effects.
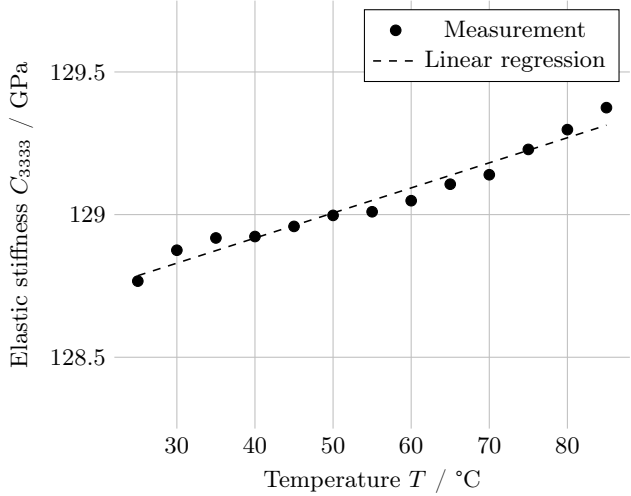
\begin{figure}
	\centering
	\begin{tikzpicture}
	\footnotesize
	\begin{axis}[
		xlabel={Temperature \(T\) / \textdegree C},
		xmajorgrids, xmin=22, xmax=88,
		ylabel={Elastic stiffness \(C_{3333}\) / GPa},
		ymajorgrids, ymin=128.25, ymax=129.75,
		axis line style={draw=none}, tick style={draw=none},
	]
		\addplot [mark size=2, only marks]
		table {%
		25 128.766841163996
		30 128.875167538772
		35 128.91812404831
		40 128.923183144298
		45 128.958666380864
		50 128.997249472252
		55 129.009953788451
		60 129.048763314288
		65 129.106575391701
		70 129.139757448561
		75 129.228507875337
		80 129.297836262438
		85 129.374893410403
		};
		\addlegendentry{Measurement}
		\addplot [semithick, dashed]
		table {%
		25 128.785910913555
		85 129.313399738702
		};
		\addlegendentry{Linear regression}
	\end{axis}
	\end{tikzpicture}
	\caption{Measurement result for the temperature dependence of the elastic parameter \(C_{3333}\) of a piezoelectric ceramic, showing a near-linear relationship~\cite{Friesen2023}.}\label{fig:elastic_measurement}
\end{figure}
Now the main novelty to be considered in the present study stems from the observation that in some materials relevant to
applications, temperature dependencies of the elastic parameters are not negligible.
As an example, a parameter of the elasticity $C$ of a piezoelectric material is shown in Fig.~\ref{fig:elastic_measurement}, indicating a near-linear increase over temperature within the considered range~\cite{Friesen2023}.
Thus led to considering situations when
\bas
	C=C(\Theta),
\eas
we note that allowing for such types of dependencies further increases the complexity with respect to modelling the thermoelastic behavior of a component in general, and to thermal losses generated by acoustic waves in particular.
In fact, it is to be expected that the overall thermal stability of a system may depend on a specific type of temperature dependence.\abs
{\bf Specifying a class of initial-boundary value problems. Notation and main results.} \quad
On supplementing (\ref{eq:visc_wave}) and (\ref{eq:thermo_coupled}) with prototypically simple
boundary and initial conditions we arrive at the problem
\be{0}
	\lball
	u_{tt} = \div (\gamma(\Theta) : \nas u_t) + a\div (\gamma(\Theta) : \nas u),
	\qquad & x\in\Om, \ t>0, \\[1mm]
	\Theta_t = D\Del\Theta + \lan \Gamma(\Theta) : \nas u_t , \nas u_t \ran,
	\qquad & x\in\Om, \ t>0, \\[1mm]
	u=0, \quad \Theta=0,
	\qquad & x\in\pO, \ t>0, \\[1mm]
	u(x,0)=u_0(x), \quad u_t(x,0)=u_{0t}(x), \quad \Theta(x,0)=\Theta_0(x),
	\qquad & x\in\Om,
	\ear
\ee
as the precise mathematical object to be subsequently considered;
in comparison with the above, to prepare convenient notation throughout our analysis we have thus set
$a=\frac{1}{\tau}$,
$D=\frac{\lambda}{c\rho}$,
$\gamma(\Theta)=\frac{\tau}{\rho} C(\Theta)$
and $\Gamma(\Theta)=\frac{\tau}{c\rho}C(\Theta)$,
although we will not rely on any relation between $\gamma$ and $\Gamma$.
We note that (\ref{0}) assumes temperature dependencies to be limited to the elasticity tensor, hence neglecting
possible variations of the parameter $a$ with respect to $\Theta$.
Apart from that, we underline that effects of thermoelasticity are not regarded within the scope of this study;
in fact, including such mechanisms is well-known to bring about substantial challenges for the mathematical analysis
already in cases when $\gamma$ does not depend on $\Theta$ (\cite{roubicek}, \cite{blanchard_guibe}).
After all, at least in simple near-linear situations in which mechanical processes are purely harmonic, such mechanisms
would lead to additional contributions which with respect to heat generation would involve production rates
with vanishing temporal averages.
In this sense, their neglection might be expected to be of minor impact in comparison to possible effects exerted by
temperature dependencies in the key system constituent $\gamma$.\abs
%
%
%
%
In order to formulate our results, and for further reference below, let us comment on the notation used throughout sequel.
For matrices $Y =(Y_{ij})_{i,j\in\{1,...,n\}}\in \R^{n\times n}\in \R^{n\times n}$ and tensors
$\beta=(\beta_{ijkl})_{i,j,k,l\in\{1,...,n\}} \in \R^{n\times n\times n\times n}$, we write
$\lan X,Y \ran := \sum_{i,j=1}^n X_{ij} Y_{ij}\in\R$
and $X^t := (X^t)_{ij}$ with $(X^t)_{ij}:=X_{ji}$ for $i,j\in\{1,...,n\}$,
and define the matrix $\beta:X\in \R^{n\times n}$ by letting
$(\beta:X)_{ij} := \sum_{k,l=1}^n \beta_{ijkl} X_{kl}$ for $i,j\in\{1,...,n\}$.
Moreover, given vectors $w=(w_1,...,w_n)\in \R^n$ and $z=(z_1,...,z_n)\in\R^n$
we introduce a matrix $w\mult z\in \R^{n\times n}$
by writing $(w\mult z)_{ij}:=w_i z_j$ for $i,j\in\{1,...,n\}$.
Finally, for vector functions $\psi=(\psi_1,...,\psi_n)\in W^{1,1}(\Om;\R^n)$ we let
$\na \psi := (\pa_j \psi_i)_{i,j=1,...,n}$ and $\nas \psi := \frac{1}{2}(\pa_j \psi_i + \pa_i \psi_j)_{i,j=1,...,n}$
denote its Jacobian and the associated symmetrized gradient, respectively.\abs
Now in an attempt to undertake a first step toward an understanding of possible influences that temperature dependencies
of the above type may have, this manuscript focuses on the issue of global solvability in (\ref{0}).
We again recall here that in the case when both $\gamma$ and $\Gamma$ are temperature-independent, then previous
literature asserts global existence of so-called weak-renormalized solutions actually even in slightly more complex
variants of (\ref{0}) in which, inter alia, some suitably mild thermoelastic effects can be admitted (\cite{blanchard_guibe};
cf.~also \cite{roubicek} for a close relative addressing Neumann bounary conditions for $\Theta$);
we emphasize, however, that precedent studies in essential parts rely on a favorable energy structure which in the
corresponding version of the simple system (\ref{0}) with $\Theta$-independent $\gamma$ is formally expressed in the identity
\be{energy}
	\frac{d}{dt} \bigg\{ \frac{1}{2} \io |u_t|^2 + \frac{a}{2} \io \lan \gamma : \nas u, \nas u\ran \bigg\}
	= - \io \lan \gamma : \nas u_t, \nas u_t\ran.
\ee
Our main results now make sure that although this structure apparently breaks down when $\gamma=\gamma(\Theta)$,
within a suitably generalized concept of solvability, to be discussed in more detail in Section \ref{sect2},
some global solution to (\ref{0}) can be found under assumptions mild enough so as to be satisfied by any sufficiently smooth
and bounded ingredients $\gamma$ and $\Gamma$:
\begin{theo}\label{theo13}
  Let $n\ge 1$ and $\Om\subset\R^n$ be a bounded domain with smooth boundary, let $a>0$ and $D>0$, and suppose
  that
  \be{gamma_Gamma_reg}
	\lbal
	\gamma=(\gamma_{ijkl})_{i,j,k,l\in\{1,...,n\}} \in
	L^\infty([0,\infty);\R^{n\times n\times n\times n}) \cap
	C^2([0,\infty);\R^{n\times n\times n\times n})
	\qquad \mbox{and} \\[1mm]
	\Gamma=(\Gamma_{ijkl})_{i,j,k,l\in\{1,...,n\}} \in
	L^\infty([0,\infty);\R^{n\times n\times n\times n}) \cap
	C^2([0,\infty);\R^{n\times n\times n\times n})
	\ear
  \ee
  are such that
  \be{gamma_symm}
	\gamma_{ijkl}(\xi)=\gamma_{klij}(\xi)
	\qquad \mbox{for all $\xi\ge 0$ and } (i,j,k,l)\in\{1,...,n\}^4
  \ee
  and
  \be{gamma2}
	\gamma_{ijkl}(\xi)=\gamma_{jikl}(\xi)
	\qquad \mbox{for all $\xi\ge 0$ and } (i,j,k,l)\in\{1,...,n\}^4,
  \ee
  that
  \be{Gamma_symm}
	\Gamma_{ijkl}(\xi)=\Gamma_{klij}(\xi)
	\qquad \mbox{for all $\xi\ge 0$ and } (i,j,k,l)\in\{1,...,n\}^4,
  \ee
  and that
  \be{gamma_pos}
	\lan \gamma(\xi):X,X\ran \ge K_\gamma |X|^2
	\qquad \mbox{for all $\xi\ge 0$ and } X\in \R^{n\times n}
  \ee
  as well as
  \be{Gamma_pos}
	\lan \Gamma(\xi):X,X\ran \ge K_\Gamma |X|^2
	\qquad \mbox{for all $\xi\ge 0$ and } X\in \R^{n\times n}
  \ee
  with some $K_\gamma>0$ and $K_\Gamma>0$.
  Then whenever
  \be{init}
	\lbal
	u_0\in W_0^{1,2}(\Om;\R^n), \\[1mm]
	u_{0t}\in L^2(\Om;\R^n) \quad \mbox{and} \\[1mm]
	\Theta_0\in L^1(\Om;\R)
	\mbox{ is nonnegative,}
	\ear
  \ee
  there exist functions
  \be{reg}
	\lbal
	u\in L^\infty_{loc}([0,\infty);W_0^{1,2}(\Om;\R^n)) \qquad \mbox{and} \\[1mm]
	\Theta\in L^\infty_{loc}([0,\infty);L^1(\Om;\R))
		\cap \bigcap_{q\in [1,\frac{n+2}{n})} L^q_{loc}(\bom\times [0,\infty);\R)
		\cap \bigcap_{r\in [1,\frac{n+2}{n+1})} L^r_{loc}([0,\infty);W_0^{1,r}(\Om;\R))
	\ear
  \ee
  such that
  \be{reg2}
	u_t\in L^\infty_{loc}([0,\infty);L^2(\Om;\R^n))
		\cap L^2_{loc}([0,\infty);W_0^{1,2}(\Om;\R^n)),
  \ee
  and that $(u,\Theta)$ forms a global generalized solution of (\ref{0}) in the sense of Definition \ref{dw} below.
\end{theo}
We note that here the assumption in (\ref{gamma2}), as already required in previous literature
(\cite{roubicek}, \cite{mielke_roubicek}, \cite{blanchard_guibe}) encodes the requirement from (\ref{stressstrainsymmetric})
that $\gamma$ turns symmetric strains to symmetric stresses.\abs
{\bf Main ideas.} \quad
Due to the absence of an energy structure of the form in (\ref{energy}),
a major challenge for any analysis of (\ref{0}) seems to be linked to the essentially quadratic type of
coupling therein: Already in the simple
particular case when $\Gamma$ coincides with the identity mapping on $\R^{n\times n}$, heat production occurs at a rate which
is bounded from below by $|\nas u_t|^2$.
In the present setting,
{\em a priori} knowledge for this quantity which is available without requirements on information about $\Theta$ seems
to reduce to estimates of the form
\bas
	\int_0^T \io \lan \gamma(\Theta):\nas (u_t+au),\nas (u_t+au)\ran \le C(T)
\eas
which reflect the dissipative action in the first equation from (\ref{0}) in the framework of some associated
zero-order energy inequality;
even in the favorable setting of bounded and uniformly positive definite $\gamma$ addressed
in Theorem \ref{theo13}, this seems to merely yield $L^2$ bounds for $\na u_t$ -- and hence, equivalently, for $\nas u_t$ --
in the sense of implying that
\bas
	\bigg( \int_0^T \io |\na u_t|^2 + \bigg) \, \int_0^T \io |\nas u_t|^2 \le C(T)
\eas
for $T>0$ (cf.~Lemma \ref{lem5}).
Accordingly, the heat source can apparently be controlled in the non-reflexive space $L^1(\Om\times (0,T))$ for $T>0$ only,
which seems to go along with a lack of knowledge on suitable compactness properties in any meaningful approximation scheme,
and hence seems to mark a crucial difference to the situation of $\Theta$-independent $\gamma$ addressed in
\cite{blanchard_guibe} and \cite{roubicek}, for instance. \abs
To overcome related difficulties,
the existence theory to be developed below resorts to a notion of generalized solvability
that substantially deviates from those introduced in some precedent studies in which certain
renormalized solution concepts for the corresponding temperature distributions were introduced, but in which
standard weak solvability is considered with respect to the displacement variable
(\cite{blanchard_guibe}, \cite{chelminski_owczarek}, \cite{chelminski_owczarek_I}).
Specifically,
with regard to the crucial solution component $\Theta$ our concept will require validity of two {\em inequalities},
instead of fulfillment of one identity such as in standard weak solution concepts.
In contrast to some cases of rather far relatives in the recent literature on fully parabolic problems
(\cite{lankeit_win_NoDEA}, \cite{win_M3AS_scrounge}),
the approach pursued here in this regard will need to appropriately cope with the wave type structure of the first equation
in (\ref{0}), and with thus fairly restricted options to make use of dissipative features.
Specifically, the core of our analysis, to be foreshadowed in Defnition \ref{dw} and executed in Lemma \ref{lem12},
will examine quantities of the form
\be{en}
	\F := \frac{1}{2} |u_t|^2 + \frac{\kappa}{2} |\na u|^2 + \lam \Theta,
	\qquad \kappa>0, \ \lam>0,
\ee
which at a pointwise level
couple the temperature field to the constituents $|u_t|^2$ and $|\na u|^2$ of the fundamental energy structure
in the wave part of (\ref{0}).
A key observation will reveal that when here the free parameters $\kappa$ and $\lam$ are chosen suitably, certain integrated
versions of $\F$ will indeed satisfy a one-sided inequality which can be viewed as providing an upper bound for
$\pa_t \F$ that is optimal in the sense of being satisfied as an identity along smooth trajectories
(see (\ref{wF}), Proposition \ref{prop_dw} and (\ref{12.99})).\abs
At a technical level, key parts of these considerations will rely on appropriate exploitation of weak lower semicontinuity
of norms in $L^2$ spaces, to be adapted to settings in which differences between such expressions and integrals
involving terms of the form
\be{quad}
	\lan \beta(z):\nas w,\nas w\ran
\ee
occur,
comparable to those in (\ref{gamma_symm}) and (\ref{gamma_pos}) (Section \ref{sect_lsc}).
Apart from that, due to limited information on regularity of $\frac{\pa\Theta}{\pa\nu}$ on $\pO$
(cf.~Lemma \ref{lem4} and Lemma \ref{lem6}), spatial localization of the arguments related to $\F$ seems in order (see (\ref{F})).
\mysection{A concept of generalized solvability. Approximate solutions}\label{sect2}
The following describes the notion of generalized solvability that will form the target object of our subsequent considerations.
While essentially standard requirements on natural weak solvability with regard to the first equation in (\ref{0}) are imposed,
with respect to the component $\Theta$ an associated one-sided inequality, (\ref{wt}),
is combined with the localized energy dissipation feature (\ref{wF}):
\begin{defi}\label{dw}
  Let $\gamma\in L^\infty([0,\infty);\R^{n\times n\times n\times n})$ and
  $\Gamma\in L^\infty([0,\infty);\R^{n\times n\times n\times n})$,
  let $a>0$ and $D>0$, and suppose that $u_0\in W_0^{1,2}(\Om;\R^n)$, $u_{0t}\in L^2(\Om;\R^n)$ and $\Theta_0\in L^1(\Om;\R)$.
  Then a pair $(u,\Theta)$ of functions
  \be{w1}
	\lbal
	u\in L^2_{loc}([0,\infty);W_0^{1,2}(\Om;\R^n)) \qquad \mbox{and} \\[1mm]
	\Theta\in L^1_{loc}([0,\infty);W_0^{1,1}(\Om;\R))
	\ear
  \ee
  will be called a {\em global generalized solution} of (\ref{0}) if
  \be{w2}
	u_t\in L^2_{loc}([0,\infty);W_0^{1,2}(\Om;\R^n))
  \ee
  and $\Theta\ge 0$ a.e.~in $\Om\times (0,\infty)$, if
  \bea{wu}
	& & \hs{-30mm}
	\int_0^\infty \io u\cdot\vp_{tt}
	+ \io u_0\cdot \vp_t(\cdot,0)
	- \io u_{0t} \cdot \vp(\cdot,0) \nn\\
	&=& - \int_0^\infty \io \lan \gamma(\Theta):\nas u_t,\na\vp \ran
	- a \int_0^\infty \io \lan \gamma(\Theta):\nas u, \na\vp\ran
  \eea
  for all $\vp\in C_0^\infty(\Om\times [0,\infty);\R^n)$, if
  \be{wt}
	- \int_0^\infty \io \Theta\wh{\vp}_t - \io \Theta_0 \wh{\vp}(\cdot,0)
	\ge D \int_0^\infty \io \Theta \Del\wh{\vp}
	+ \int_0^\infty \io \lan \Gamma(\Theta):\nas u_t,\nas u_t \ran \wh{\vp}
  \ee
  for each nonnegative $\wh{\vp}\in C_0^\infty(\Om\times [0,\infty);\R)$, and if there exist $\kappa>0, \lam>0$ and $\mu>0$
  with the property that for any nonnegative $\psi\in C_0^\infty(\Om;\R)$ and arbitrary nonincreasing
  $\zeta\in C_0^\infty([0,\infty);\R)$, the inequality
  \be{wF}
	\int_0^\infty \io \D^{(\kappa,\lam,\mu,\zeta,\psi)}
	\le \zeta(0) \io \F_0^{(\kappa,\lam,\psi)}
	+ \int_0^\infty \io \Rest^{(\lam,\mu,\zeta,\psi)}
  \ee
  holds, where
  \be{F0}
	\F_0^{(\kappa,\lam,\psi)}
	\equiv \F_0^{(\kappa,\lam,\psi)}(x)
	:= \Big(\frac{1}{2} |u_{0t}|^2 + \frac{\kappa}{2} |\na u_0|^2 + \lam \Theta_0\Big) \psi,
	\qquad x\in\Om,
  \ee
  and
  \bea{D}
	\D^{(\kappa,\lam,\mu,\zeta,\psi)}
	&\equiv& \D^{(\kappa,\lam,\mu,\zeta,\psi)} (x,t) \nn\\
	&:=&
	\Big\{
	\lan \gamma(\Theta):\nas u_t, \nas u_t \ran
	+ a \lan \gamma(\Theta):\nas u, \nas u_t \ran \nn\\
	& & \hs{4mm}
	- \kappa \lan \na u,\na u_t \ran
	- \lam \lan \Gamma(\Theta):\nas u_t, \nas u_t \ran
	\Big\} \, \zeta(t) e^{-\mu t} \psi \nn\\
	& & + \F^{(\kappa,\lam,\psi)} \big(\mu\zeta(t)-\zeta_t(t)\big) e^{-\mu t},
	\qquad x\in\Om, \ t>0,
  \eea
  with
  \be{F}
	\F^{(\kappa,\lam,\psi)}
	\equiv \F^{(\kappa,\lam,\psi)}(x,t)
	:=\Big(\frac{1}{2} |u_t|^2 + \frac{\kappa}{2} |\na u|^2 + \lam\Theta\Big) \psi,
	\qquad x\in\Om, \ t>0,
  \ee
  and where
  \bea{R}
	\Rest^{(\lam,\mu,\zeta,\psi)}
	&\equiv& \Rest^{(\lam,\mu,\zeta,\psi)}(x,t) \nn\\
	&:=& \Big\{
	\lan \gamma(\Theta):\nas u_t, u_t\mult \na \psi \ran
	+ a \lan \gamma(\Theta):\nas u, u_t\mult \na \psi \ran
	+ \lam D \Theta \Del\psi
	\Big\} \, \zeta(t) e^{-\mu t}, \nn\\[0mm]
	& & \hs{90mm}
	\qquad x\in\Om, \ t>0.
  \eea
\end{defi}
Consistency of this concept with that of classical solvability is underlined by the following observation.
%
%
\begin{prop}\label{prop_dw}
  Let $\gamma$ and $\Gamma$ belong to
  $L^\infty([0,\infty);\R^{n\times n\times n\times n}) \cap C^1([0,\infty);\R^{n\times n\times n\times n})$,
  let $a>0$ and $D>0$,
  and let $u_0\in C^0(\bom;\R^n), u_{0t} \in C^0(\bom;\R^n), \Theta_0\in C^0(\bom;\R)$ as well as
  \bas
	\lbal
	u \in C^0(\bom\times [0,\infty);\R^n)
		\cap C^{2,1}(\Om\times (0,\infty);\R^n) \qquad \mbox{and} \\[1mm]
	\Theta\in C^0(\bom\times [0,\infty);\R)
		\cap C^{2,1}(\Om\times (0,\infty);\R)
	\ear
  \eas
  be such that
  \bas
	u_t \in C^0(\bom\times [0,\infty);\R^n)
		\cap C^{2,1}(\Om\times (0,\infty);\R^n)
  \eas
  \begin{align}
  &t\mapsto \io |∇u(\cdot,t)|^2 \quad\text{ is continuous, }\nn\\
  & u_{tt}\in \Lloc1(\Omega\times[0,∞)),\quad Θ_t, ΔΘ\in \Lloc1(\Omega\times[0,∞))\label{eq:integrability}
  \end{align}
  and that $(u,\Theta)$ forms a global generalized solution of (\ref{0}) in the sense of Definition \ref{dw}.
  Then (\ref{0}) is actually satisfied in the classical sense.
\end{prop}
\proof
That $u=0$ and $\Theta=0$ on $\bom\times(0,\infty)$ immediately results from the functions belonging to $\Lloc2([0,∞);W_0^{1,2}(\Om; ℝ^n))$ and $\Lloc1([0,∞);W_0^{1,1}(\Om;ℝ^n)$, respectively, and their continuity.

If we insert arbitrary $φ\in C_0^{∞}(\Om\times(0,∞))$ in \eqref{wu} and integrate by parts twice with respect to time and once to space (which is obviously possible, because all of the integrands $φu_{tt}$, $φ_tu_t$, $φ_{tt} u$, $\lan \gamma(\Theta):\nas u_t,\na\vp \ran + a \lan \gamma(\Theta):\nas u, \na\vp\ran $, $\div(γ(Θ):\nas (u_t+au)) φ$ are continuous on $\supp φ$), from the fundamental lemma of the calculus of variations (applicable, since $u_{tt}$ and $\div(γ(Θ):\nas (u_t+au)$ belong to $\Lloc1(\Om\times(0,∞))$), we obtain that the first equation of \eqref{0} is satisfied at every point in $\Omega\times(0,∞)$. For $ψ\in C_0^\infty(\Om)$ and $ε>0$ inserting $φ(x,t)=ψ(x)(\f{t}{ε}-\f{2t^2}{ε^2}+\f{t^3}{ε^3})χ_{[0,ε]}(t)$ or $φ(x,t)=ψ(x)(1-\f{2t^2}{ε^2}+\f{t^4}{ε^4})χ_{[0,ε]}(t)$ in \eqref{wu}
and taking the limit $ε\searrow 0$ (using continuity of $u$ and $u_t$ at $t=0$) shows that $u$ fulfils the initial conditions.

From \eqref{eq:integrability} and an analogous treatment of \eqref{wt} we find that $Θ_t\ge DΔΘ + \lan \Gamma(\Theta):\nas u_t,\nas u_t \ran$ in $\Om\times(0,∞)$ and $Θ(0)\ge Θ_0$ in $\Om$.

From a choice of $ζ(t):= (1-\f{t}{ε})_+$ in \eqref{wF} due to continuity of $t\mapsto \io \F^{(\kappa,\lam,\psi)}(\cdot,t)$,
\[
\io\Big(\frac{1}{2} |u_t(\cdot,0)|^2 + \frac{\kappa}{2} |\na u(\cdot,0)|^2 + \lam\Theta(\cdot,0)\Big) \psi \le \io \Big(\frac{1}{2} |u_{0t}|^2 + \frac{\kappa}{2} |\na u_0|^2 + \lam \Theta_0\Big) \psi
\]
for each nonnegative $ψ\in C_0^\infty(\Om)$ and hence $Θ(0)\le Θ_0$ in $\Om$.

Accordingly, for any $ψ\in C_0^\infty(\Om)$ and $ζ\in C_0^\infty([0,∞))$,
\newcommand{\kl}[1]{\left(#1\right)}
\begin{align*}
 \int_0^\infty\io &\kl{\f12|u_t|^2 + \f{κ}2 |∇u|^2 +λΘ} ψ(-ζe^{-μt})_t =\\&\int_0^\infty \io \kl{u_tu_{tt}+κ\lan ∇u,\nas u_t\ran + λΘ_t} ψζe^{-μt} + ζ(0) \io \kl{|u_{0t}|^2 +\f{κ}2 |∇u_0|^2 +λΘ_0} ψ
\end{align*}
since $u_{tt}, Θ_t\in \Lloc1(\Omega\times[0,∞))$.
If we insert the equality for $u$ and integrate by parts with respect to space, many terms in \eqref{wF} are cancelled and we see that
\[
 \int_0^{∞}\io
	\kl{- \lam \lan \Gamma(\Theta):\nas u_t, \nas u_t \ran + λΘ_t} \, \zeta(t) e^{-\mu t} \psi
\le
 \int_0^\infty\io
	 \lam D \Theta \Del\psi \zeta(t) e^{-\mu t},
\]
which finally due to $ΔΘ\in \Lloc1(\Om\times[0,∞))$ shows the remaining inequality $Θ_t\le DΔΘ + \lan \Gamma(\Theta):\nas u_t,\nas u_t \ran$.

%

\qed
%
%
%
%
%
%
%
%
%
%
%
%
%
%
%
%
%
%
%
Throughout the sequel, we shall consider $\Om, a,D,\gamma$ and $\Gamma$ as well as $u_0,u_{0t}$ and $\Theta_0$ to be fixed and such that the assumptions of Theorem \ref{theo13} are met.
We can then pick $(u_{0\eps})_{\eps\in (0,1)} \subset C_0^\infty(\Om;\R^n)$,
$(u_{0t\eps})_{\eps\in (0,1)} \subset C_0^\infty(\Om;\R^n)$ and
$(\Theta_{0\eps})_{\eps\in (0,1)} \subset C_0^\infty(\Om;\R)$
in such a way that $\Theta_{0\eps} \ge 0$ in $\Om$ for all $\eps\in (0,1)$, and that
\be{ie}
	u_{0\eps} \to u_0
	\quad \mbox{in } W^{1,2}(\Om;\R^n),
	\quad
	u_{0t\eps} \to u_{0t}
	\quad \mbox{in } L^2(\Om;\R^n)
	\quad \mbox{and} \quad
	\Theta_{0\eps} \to \Theta_0
	\quad \mbox{in } L^1(\Om;\R)
	\qquad \mbox{as } \eps\searrow 0.
\ee
In line with standard theory of local solvability in parabolic systems (\cite{amann}),
this particularly ensures that
for each $\eps\in (0,1)$, the regularized variant of (\ref{0}) given by
\be{0eps}
	\lball
	v_{\eps t} = - \eps \Del^2 \veps + \div (\gamma(\Teps) : \nas \veps) + a\div (\gamma(\Teps) : \nas \ueps),
	\qquad & x\in\Om, \ t>0, \\[1mm]
	u_{\eps t} = \eps \Del \ueps + \veps,
	\qquad & x\in\Om, \ t>0, \\[1mm]
	\Theta_{\eps t} = D\Del\Teps + \lan \Gamma(\Teps) : \nas \veps, \nas \veps \ran,
	\qquad & x\in\Om, \ t>0, \\[1mm]
	\veps=0, \quad Δ\veps=0,\quad \ueps=0, \quad \Teps=0,
	\qquad & x\in\pO, \ t>0, \\[1mm]
	\veps(x,0)=u_{0t\eps}(x), \quad \ueps(x,0)=u_{0\eps}(x), \quad \Teps(x,0)=\Theta_{0\eps}(x),
	\qquad & x\in\Om,
	\ear
\ee
admits a local-in-time classical solution in the following sense:
\begin{lem}\label{lem_loc}
  Let $\eps\in (0,1)$. Then there exist $\tme\in (0,\infty]$ as well as functions
  \bas
	\lbal
	\veps\in C^0(\bom\times [0,\tme);\R^n)
		\cap C^{4,1}(\bom\times (0,\tme);\R^n), \\[1mm]
	\ueps\in C^0(\bom\times [0,\tme);\R^n)
		\cap C^{2,1}(\bom\times (0,\tme);\R^n)
		\cap C^0([0,\tme);W_0^{1,2}(\Om;\R^n)) \qquad \mbox{and} \\[1mm]
	\Teps \in C^0(\bom\times [0,\tme);\R)
		\cap C^{2,1}(\bom\times (0,\tme);\R)
	\ear
  \eas
  such that $\Teps\ge 0$ in $\bom\times [0,\tme)$, that $(\veps,\ueps,\Teps)$ solves (\ref{0eps})
  classically in $\Om\times (0,\tme)$, and that
  \bea{ext}
	& & \hs{-20mm}
	\mbox{if $\tme<\infty$, \quad then for all $\eta>0$,} \nn\\
	& & \hs{-10mm}
	\limsup_{t\nearrow \tme} \Big\{
	\|\veps(\cdot,t)\|_{W^{2+2\eta,\infty}(\Om)}
	+ \|\ueps(\cdot,t)\|_{W^{1+\eta,\infty}(\Om)}
	+ \|\Teps(\cdot,t)\|_{W^{1+\eta,\infty}(\Om)}
	\Big\} = \infty
  \eea
\end{lem}
\proof
For any $p>n$, we let $E_1=\set{(v,u,Θ)\in W^{4,p}(\Om)\times W^{2,p}(\Om)\times W^{2,p}(\Om)\mid v\bdry=0, Δv\bdry=0, u\bdry=0, Θ\bdry=0}$ and $E_0=(\Lom{p})^3$ and denote by $E_{γ}=(E_0,E_1)_{γ}$ the interpolation space of order $γ\in(0,1)$, let $f\equiv 0$ and
\[A(v,u,Θ)=\begin{pmatrix}εΔ^2-\div(γ(Θ):\nas \cdot)&-a\div(γ(Θ)\nas\cdot)&\\&-εΔ&\\-\lan Γ(Θ):\nas v:\nas \cdot\ran&&-DΔ\end{pmatrix}\in L(E_1,E_0), \qquad (u,v,Θ)\in E_{β}.\]
Since for $(v_1,u_1,Θ_1),(v_2,u_2,Θ_2)\in E_{β}$ and $(\tilde v,\tilde u,\tilde{Θ})\in E_1$,
\begin{align*}
 &\norm[\Lom p]{\div((γ(Θ_1)-γ(Θ_2))\nas \tilde u)}\\
\le  &\norm[\Lom p]{((γ(Θ_1)-γ(Θ_2))D^2 \tilde u)}
+ \norm[\Lom p]{(γ'(Θ_1)-γ'(Θ_2))∇Θ_1 \nas \tilde u)}
+ \norm[\Lom p]{γ'(Θ_2)
(∇Θ_1-∇Θ_2) \nas \tilde u}\\
\le &(\norm[\Lom \infty]{γ(Θ_1)-γ(Θ_2)}
+ \norm[\Lom p]{(γ'(Θ_1)-γ'(Θ_2))∇Θ_1)}
+ \norm[\Lom p]{γ'(Θ_2)
(∇Θ_1-∇Θ_2)})\norm[W^{2,p}]{\tilde u}
\end{align*}
and thus
\begin{align*}
 &\norm[L(E_1,E_0)]{\div((γ(Θ_1)-γ(Θ_2))\nas\cdot)}\\
 &\le \norm[\Lom \infty]{γ(Θ_1)-γ(Θ_2)}
+ \norm[\Lom p]{(γ'(Θ_1)-γ'(Θ_2))∇Θ_1}
+ \norm[\Lom p]{γ'(Θ_2)
(∇Θ_1-∇Θ_2)}
\end{align*}
and related estimates for the other terms
show the required Lipschitz continuity of $A\colon E_{β}\to L(E_1,E_0)$ for $β=\f12$, we may employ the general existence result of \cite[Thm.~12.1]{amann}, obtaining a solution and deriving an extensibility criterion in $E_{δ}$ for $δ>\f12$ and thus, in consequence, \eqref{ext}, from \cite[Thm.~12.5]{amann}.

\qed
\mysection{Basic testing procedures. {\em A priori} estimates for \tops{$\veps$}{v eps} and \tops{$\ueps$}{u eps}}
A first testing procedure applied to (\ref{0eps}) is designed here in such a way that not only the
derivation of spatially global estimates is prepared (Lemma \ref{lem5}), but that later on
also our analysis of (\ref{wF}) can be built on this (see Lemma \ref{lem12}).
\begin{lem}\label{lem2}
  Let $\psi\in C^2(\bom;\R)$ and $\eps\in (0,1)$. Then
  \bea{2.1}
	& & \hs{-20mm}
	\frac{1}{2} \io (|\veps|^2)_t \psi
	+ \io \lan \gamma(\Teps):\nas \veps,\nas\veps\ran \psi
	+ \eps \io |\Del\veps|^2 \psi \nn\\
	&=& - a \io \lan\gamma(\Teps):\nas\ueps,\nas\veps\ran\psi \nn\\
	& & \jl{-} \io \lan\gamma(\Teps):\nas\veps,\veps\mult\na\psi\ran
	\jl{-} a \io \lan\gamma(\Teps):\nas\ueps,\veps\mult\na\psi\ran \nn\\
	& & - 2\eps \io (\na\veps\cdot\na\psi)\cdot\Del\veps
	- \eps \io (\veps\cdot\Del\veps) \Del\psi
	\qquad \mbox{ for all } t\in (0,\tme).
  \eea
\end{lem}
\proof
  According to the first equation in (\ref{0eps}),
  \bea{2.2}
	\frac{1}{2} \io (|\veps|^2)_t \psi
	&=& \io \veps \cdot \div \big(\gamma(\Teps):\nas\veps\big) \psi
	+ a \io \veps\cdot\div \big(\gamma(\Teps):\nas\ueps\big)\psi \nn\\
	& & - \eps \io (\veps\cdot\Del^2 \veps)\psi
	\qquad \mbox{ for all } t\in (0,\tme),
  \eea
  where thanks to the boundary conditions $\veps=0$ and $\Del\veps=0$ on $\pO\times (0,\tme)$, two integrations by parts show that
  for all $t\in (0,\tme)$,
  \bea{2.3}
	- \eps \io (\veps\cdot\Del^2 \veps)\psi
	&=& - \eps \io \Del\veps\cdot\Del (\veps\cdot\psi) \nn\\
	&=& - \eps \io |\Del\veps|^2 \psi
	- 2\eps \io (\na\veps\cdot\na\psi)\cdot\Del\veps
	- \eps \io (\veps\cdot\Del\veps)\Del\psi.
  \eea
  Again since $\veps\equiv (v_{\eps 1},...,v_{\eps_n})=0$ on $\pO\times (0,\tme)$, by another integration by parts we find that
  \bea{2.4}
	\io \veps \cdot \div \big(\gamma(\Teps):\nas\veps\big) \psi
	&=& \sum_{i=1}^n \io \psi v_{\eps i} \Big( \div\big(\gamma(\Teps):\nas\veps\big) \Big)_i \nn\\
	&=& \sum_{i,j=1}^n \io \psi v_{\eps i} \pa_j \big(\gamma(\Teps):\nas\veps\big)_{ij} \nn\\
	&=& - \sum_{i,j=1}^n \io \pa_j(\psi v_{\eps i}) \big(\gamma(\Teps):\nas\veps\big)_{ij} \nn\\
	&=& - \sum_{i,j=1}^n \io \psi \pa_j v_{\eps i} \big(\gamma(\Teps):\nas\veps\big)_{ij}
	- \sum_{i,j=1}^n \io \pa_j \psi v_{\eps i} \big(\gamma(\Teps):\nas\veps\big)_{ij} \nn\\
	&=& - \io \lan\gamma(\Teps):\nas\veps,\na\veps\ran\psi
	- \io \lan\gamma(\Teps):\nas\veps,\veps\mult\na\psi\ran \nn\\
	&=& - \io \lan\gamma(\Teps):\nas\veps,\nas\veps\ran\psi
	- \io \lan\gamma(\Teps):\nas\veps,\veps\mult\na\psi\ran
  \eea
  for all $t\in (0,\tme)$,
  because our assumption (\ref{gamma2}) guarantees that
  \bas
	\lan\gamma(\Teps):\nas\veps,(\na\veps)^t\ran
	&=& \sum_{i,j,k,l=1}^n \gamma_{ijkl}(\Teps) (\pa_l v_{\eps k} + \pa_k v_{\eps l}) \pa_i v_{\eps j} \\
	&=& \sum_{i,j,k,l=1}^n \gamma_{jikl}(\Teps) (\pa_l v_{\eps k} + \pa_k v_{\eps l}) \pa_i v_{\eps j} \\
	&=& \sum_{i',j',k,l=1}^n \gamma_{i'j'kl}(\Teps) (\pa_l v_{\eps k} + \pa_k v_{\eps l}) \pa_{j'} v_{\eps i'} \\[2mm]
	&=& \lan\gamma(\Teps):\nas\veps,\na\veps \ran
  \eas
  and hence
  \bas
	\lan\gamma(\Teps):\nas\veps,\na\veps\ran
	&=& \frac{1}{2} \lan\gamma(\Teps):\nas\veps,\na\veps\ran
	+ \frac{1}{2} \lan\gamma(\Teps):\nas\veps,(\na\veps)^t\ran
	= \lan\gamma(\Teps):\nas\veps,\nas\veps\ran
  \eas
  in $\Om\times (0,\tme)$.
  As, similarly,
  \bas
	a \io \veps\cdot\div \big(\gamma(\Teps):\nas\ueps\big)\psi
	= - a \io \lan\gamma(\Teps):\nas\ueps,\nas\veps\ran
	- a \io \lan\gamma(\Teps):\nas\ueps,\veps\mult\na\psi\ran
  \eas
  for all $t\in (0,\tme)$,
  combining (\ref{2.2}) with (\ref{2.3}) and (\ref{2.4}) leads to (\ref{2.1}).
\qed
A second basic feature of (\ref{0eps}) is rather evident.
\begin{lem}\label{lem3}
  If $\psi\in C^1(\bom;\R)$ and $\eps\in (0,1)$, then
  \be{3.1}
	\frac{1}{2} \io \big( |\na\ueps|^2\big)_t \psi
	+ \eps \io |\Del\ueps|^2 \psi
	= \io \lan\na\ueps,\na\veps\ran \psi
	- \eps \io (\na\ueps\cdot\na\psi)\cdot\Del\ueps
  \ee
  for all $t\in (0,\tme)$.
\end{lem}
\proof
  From the second equation in (\ref{0eps}), we obtain the identity
  \be{3.2}
	\frac{1}{2} \io \big(|\na\ueps|^2\big)_t \psi
	= \io \lan\na\ueps,\na\veps\ran\psi
	+ \eps \io \lan \psi\na\ueps,\na\Del\ueps\ran
	\qquad \mbox{ for all } t\in (0,\tme).
  \ee
  Since the boundary condition $\ueps=(u_{\eps 1},...,u_{\eps n})=0$ on $\pO\times (0,\tme)$
  implies that also $u_{\eps t}|_{\pO\times (0,\tme)}=0$,
  and since thus $\Del\ueps=\frac{u_{\eps t}-\veps}{\eps}=0$ on $\pO\times (0,\tme)$, we may here integrate by parts without
  encountering nonzero boundary integrals, thereby confirming that
  \bas
	\eps \io \lan \psi\na\ueps,\na\Del\ueps\ran
	&=& \eps \sum_{i,j=1}^n \io \psi \pa_j u_{\eps i} \pa_j \Del u_{\eps i} \\
	&=& - \eps \sum_{i,j=1}^n \io \psi \pa_{jj} u_{\eps i} \Del u_{\eps i}
	- \eps \sum_{i,j=1}^n \io \pa_j \psi \pa_j u_{\eps i} \Del u_{\eps i} \\
	&=& - \eps \io |\Del\ueps|^2 \psi
	- \eps \io (\na\ueps\cdot\na\psi)\cdot\Del\ueps
	\qquad \mbox{ for all } t\in (0,\tme).
  \eas
  Therefore, (\ref{3.2}) is equivalent to (\ref{3.1}).
\qed
We furthermore record another simple property of solutions to (\ref{0eps}).
\begin{lem}\label{lem4}
  Whenever $\psi\in C^2(\bom;\R)$ and $\eps\in (0,1)$, we have
  \be{4.1}
	\io \Theta_{\eps t} \psi
	= D \io \Teps \Del\psi
	+ D \int_{\pO} \frac{\pa\Teps}{\pa\nu} \psi
	+ \io \lan\Gamma(\Teps):\nas\veps,\nas\veps\ran\psi
	\qquad \mbox{for all } t\in (0,\tme).
  \ee
\end{lem}
\proof
  From (\ref{0eps}) we obtain that
  \bas
	\io \Theta_{\eps t} \psi
	= D \io \Del\Teps \psi
	+ \io \lan\Gamma(\Teps):\nas\veps,\nas\veps\ran\psi
	\qquad \mbox{for all } t\in (0,\tme),
  \eas
  and since two integrations by parts relying on the identity $\Teps|_{\pO\times (0,\tme)}=0$ show that
  \bas
	D \io \Del\Teps\psi
	&=& - D \io \na\Teps\cdot\na\psi
	+ D\int_{\pO} \frac{\pa\Teps}{\pa\nu} \psi \\
	&=& D \io \Teps\Del\psi
	+ D\int_{\pO} \frac{\pa\Teps}{\pa\nu} \psi
	\qquad \mbox{for all } t\in (0,\tme),
  \eas
  this yields (\ref{4.1}).
\qed
Choosing $\psi\equiv 1$ in (\ref{2.1}) and (\ref{3.1}) lets a linear combination of the resulting identities
become an essentially straightforward energy analysis of the wave subsystem of (\ref{0eps}), leading to fairly natural
results as follows.
\begin{lem}\label{lem5}
  Let $T>0$. Then there exists $C(T)>0$ such that
  \be{5.1}
	\io |\veps(\cdot,t)|^2 \le C(T)
	\qquad \mbox{for all $t\in (0,T)\cap (0,\tme)$ and } \eps\in (0,1)
  \ee
  and
  \be{5.2}
	\io |\na\ueps(\cdot,t)|^2 \le C(T)
	\qquad \mbox{for all $t\in (0,T)\cap (0,\tme)$ and } \eps\in (0,1),
  \ee
  that
  \be{5.4}
	\int_0^t \io |\na\veps|^2 \le C(T)
	\qquad \mbox{for all $t\in (0,T)\cap (0,\tme)$ and } \eps\in (0,1)
  \ee
  and that
  \be{5.5}
	\eps \int_0^t \io |\Del\veps|^2 \le C(T)
	\qquad \mbox{for all $t\in (0,T)\cap (0,\tme)$ and } \eps\in (0,1)
  \ee
  as well as
  \be{5.6}
	\eps \int_0^t \io |\Del\ueps|^2 \le C(T)
	\qquad \mbox{for all $t\in (0,T)\cap (0,\tme)$ and } \eps\in (0,1).
  \ee
\end{lem}
\proof
  On choosing $\psi\equiv 1$ therein, from Lemma \ref{lem2} and Lemma \ref{lem3} we infer that
  \be{5.66}
	\yeps(t):=\io |\veps(\cdot,t)|^2 + \io |\na\ueps(\cdot,t)|^2,
	\qquad t\in [0,\tme), \eps\in (0,1),
  \ee
  satisfies
  \bea{5.7}
	& & \hs{-20mm}
	\yeps'(t)
	+ 2 \io \lan\gamma(\Teps):\nas\veps,\nas\veps\ran
	+ 2\eps \io |\Del\veps|^2
	+ 2\eps \io |\Del\ueps|^2 \nn\\
	&=& -2a\io \lan \gamma(\Teps):\nas\ueps,\nas\veps\ran
	+ 2\io \lan\na\ueps,\na\veps\ran
	\quad \mbox{for all $t\in (0,\tme)$ and } \eps\in (0,1).
  \eea
  Here, a combination of (\ref{gamma_pos}) with Korn's inequality (see e.g. \cite{neff_pauly_witsch,korn})
  shows that with some $c_1>0$ we have
  \bas
	2 \io \lan\gamma(\Teps):\nas\veps,\nas\veps\ran
	\ge 2K_\gamma \io |\nas\veps|^2
	\ge c_1 \io |\na\veps|^2
	\qquad \mbox{for all $t\in (0,\tme)$ and } \eps\in (0,1),
  \eas
  while the boundedness of $\gamma$ on $[0,\infty)$ ensures the existence of $c_2>0$ such that thanks to Young's inequality,
  \bas
	-2a\io \lan \gamma(\Teps):\nas\ueps,\nas\veps\ran
	&\le& c_2 \io |\na\ueps| \, |\na\veps| \\
	&\le& \frac{c_1}{4} \io |\na\veps|^2 + \frac{c_2^2}{c_1} \io |\na\ueps|^2
	\qquad \mbox{for all $t\in (0,\tme)$ and } \eps\in (0,1).
  \eas
  As moreover, again by Young's inequality,
  \be{5.8}
	2\io \lan\na\ueps,\na\veps\ran
	\le \frac{c_1}{4} \io |\na\veps|^2
	+ \frac{4}{c_1} \io |\na\ueps|^2
	\qquad \mbox{for all $t\in (0,\tme)$ and } \eps\in (0,1),
  \ee
  from (\ref{5.7}) and (\ref{5.66}) we conclude that
  \bas
	\yeps'(t) + \heps(t) \le c_3 \yeps(t)
	\qquad \mbox{for all $t\in (0,\tme)$ and } \eps\in (0,1),
  \eas
  where $c_3:=\frac{c_2^2}{c_1} + \frac{4}{c_1}$, and where
  \bas
	\heps(t)
	:= \frac{c_1}{2} \io |\na\veps(\cdot,t)|^2
	+ 2\eps \io |\Del\veps(\cdot,t)|^2
	+ 2\eps \io |\Del\ueps(\cdot,t)|^2,
	\qquad t\in (0,\tme), \eps\in (0,1).
  \eas
  By nonnegativity of $\heps$ for $\eps\in (0,1)$, using Gronwall's inequality we firstly infer from (\ref{5.8}) that
  \be{5.9}
	\yeps(t) \le c_4:=\Big\{ \sup_{\eps\in (0,1)} \yeps(0)\Big\} e^{c_3 T}
	\qquad \mbox{for all $t\in [0,T)\cap [0,\tme)$ and } \eps\in (0,1),
  \ee
  with $c_4$ being finite due to (\ref{ie}).
  This directly yields (\ref{5.1}) and (\ref{5.2}), whereas (\ref{5.4}), (\ref{5.5}) and (\ref{5.6}) result from an integration
  in (\ref{5.8}), which in view of (\ref{5.9}) shows that
  \bas
	\int_0^t \heps(s) ds
	\le \yeps(0) + c_3 \int_0^t \yeps(s) ds
	\le c_4 + c_3 c_4 T
	\qquad \mbox{for all $t\in (0,T)\cap (0,\tme)$ and } \eps\in (0,1),
  \eas
  namely.
\qed
\mysection{Global existence in the regularized problems}
For fixed $\eps\in (0,1)$, due to the fourth-order parabolic regularization in the first equation from (\ref{0eps})
the moderate information on $L^2$-boundedness of $\veps$ in (\ref{5.1}) is already sufficient to ensure bounds for this
quantity with respect to the norm in any of the spaces $W^{s,p}(\Om;\R^n)$ with $s<3$ and $p<\infty$; indeed:
\begin{lem}\label{lem55}
  Let $p\ge 2$, and let $A$ denote the realization of $\Del^2$ under the boundary conditions $(\cdot)|_{\pO}=0$ and
  $\Del(\cdot)|_{\pO}=0$ in $L^p(\Om;\R^n)$, with domain given by
  $D(A):=\{\psi\in W^{4,p}(\Om;\R^n) \ | \ \psi=\Del\psi=0 \mbox{ on } \pO \}$.
  Then whenever $\eps\in (0,1)$ is such that $\tme<\infty$, for each $\al\in (\frac{1}{4},\frac{3}{4})$ there exists
  $C(\al,\eps)>0$ such that the corresponding fractional power satisfies
  \be{55.1}
	\|A^\al \veps(\cdot,t)\|_{L^p(\Om)} \le C(\al,\eps)
	\qquad \mbox{for all $t\in (\frac{1}{2}\tme,\tme)$.}
  \ee
\end{lem}
\proof
  We first note that according to known smoothing properties of the analytic semigroup $(e^{-t\eps A})_{t\ge 0}$
  (\cite{friedman_book})
  and a duality argument in the flavor of \cite[Lemma 2.1]{horstmann_win},
  one can find $\al'=\al'(\al)\in (\frac{1}{4},\frac{3}{4})$, $c_1=c_1(\al,\eps)>0$ and $c_2=c_2(\al,\eps)>0$ such that
  \bas
	\|A^\al e^{-t\eps A} \div \psi \|_{L^p(\Om)}
	\le c_1 t^{-\al'-\frac{1}{4}} \|\psi\|_{L^p(\Om)}
	\qquad \mbox{for all $t>0$ and each $\psi\in C^1(\bom;\R^n)$ fulfilling $\psi|_{\pO}=0$,}
  \eas
  and that
  \bas
	\|A^\al e^{-t\eps A} \psi\|_{L^p(\Om)} \le c_2 \|\psi\|_{W^{4,p}(\Om)}
	\qquad \mbox{for all $t>0$ and any } \psi\in D(A).
  \eas
  Then writing $t_0:=\frac{1}{2} \tme$ and using a Duhamel representation associated with the first equation in (\ref{0eps}),
  we see that
  for all $t\in (t_0,\tme)$,
  \bas
	\|A^\al \veps(\cdot,t)\|_{L^p(\Om)}
	&=& \bigg\| A^\al e^{-(t-t_0)\eps A} \veps(\cdot,t_0)
	+ \int_{t_0}^t A^\al e^{-(t-s)\eps A} \div \big\{ \gamma(\Teps):\nas (\veps+a\ueps) \big\}(\cdot,s) ds \bigg\|_{L^p(\Om)} \\
	&\le& \|A^\al e^{-(t-t_0)\eps A} \veps(\cdot,t_0)\|_{L^p(\Om)} \\
	& & + c_1 \int_{t_0}^t (t-s)^{-\al'-\frac{1}{4}} \big\| \gamma(\Teps):\nas(\veps+a\ueps)\big\|_{L^p(\Om)} (\cdot,s) ds,
  \eas
  so that since $\veps(\cdot,t_0)\in D(A)$ by Lemma \ref{lem_loc}, and since $\gamma$ is bounded on $[0,\infty)$, with some
  $c_3=c_3(\al,\eps)>0$ we have
  \be{55.2}
	\|A^\al \veps(\cdot,t)\|_{L^p(\Om)}
	\le c_3 + c_3 \int_{t_0}^t (t-s)^{-\al'-\frac{1}{4}}
	\big\{ \|\veps(\cdot,s)\|_{W^{1,p}(\Om)} + \|\ueps(\cdot,s)\|_{W^{1,p}(\Om)} \big\} ds
  \ee
  for all $t\in (t_0,\tme)$.
  Since a standard interpolation property involving fractional powers of $A$ (\cite{friedman_book})  provides $\vt_1=\vt(\al)\in (0,1)$ and $c_4=c_4(\al)>0$ fulfilling
  \bas
	\|\veps(\cdot,s)\|_{W^{1,p}(\Om)}
	\le c_4 \|A^\al \veps(\cdot,s)\|_{L^p(\Om)}^{\vt_1} \|\veps(\cdot,s)\|_{L^p(\Om)}^{1-\vt_1}
	\qquad \mbox{for all } s\in (t_0,\tme),
  \eas
  and since the Gagliardo-Nirenberg inequality together with (\ref{5.1}) shows that with some $\vt_2\in (0,1)$, $c_5>0$ and $c_6>0$
  we have
  \bas
	\|\veps(\cdot,s)\|_{L^p(\Om)}
	\le c_5 \|\veps(\cdot,s)\|_{W^{1,p}(\Om)}^{\vt_2} \|\veps(\cdot,s)\|_{L^2(\Om)}^{1-\vt_2}
	\le c_6 \|\veps(\cdot,s)\|_{W^{1,p}(\Om)}^{\vt_2}
	\qquad \mbox{for all } s\in (t_0,\tme),
  \eas
  we readily infer the existence of $\vt_3=\vt_3(\al)\in (0,1)$ and $c_7=c_t(\al,\eps)>0$ such that if for $T\in (t_0,\tme)$ we let
  \bas
	M(T):=\sup_{t\in (t_0,T)} \|A^\al \veps(\cdot,t)\|_{L^p(\Om)},
  \eas
  then
  \be{55.3}
	\|\veps(\cdot,s)\|_{W^{1,p}(\Om)}
	\le c_7 \|A^\al\veps(\cdot,s)\|_{L^p(\Om)}^{\vt_3}
	\le c_7 M^{\vt_3}(T)
	\qquad \mbox{for all $s\in (t_0,T)$ and } T\in (t_0,\tme).
  \ee
  As $\al'+\frac{1}{4}<1$, (\ref{55.2}) therefore implies the existence of $c_8=c_8(\al,\eps)>0$ such that
  \bea{55.4}
	\|A^\al \veps(\cdot,t)\|_{L^p(\Om)}
	&\le& c_3 + c_3 c_7 M^{\vt_3}(T) \int_{t_0}^t (t-s)^{-\al'-\frac{1}{4}} ds
	+ c_3 \int_{t_0}^t (t-s)^{-\al'-\frac{1}{4}} \|\ueps(\cdot,s)\|_{W^{1,p}(\Om)} ds \nn\\
	&\le& c_8 (1+M^{\vt_3}(T))
	+ c_3 \int_{t_0}^t (t-s)^{-\al'-\frac{1}{4}} \|\ueps(\cdot,s)\|_{W^{1,p}(\Om)} ds
  \eea
  for all $t\in (t_0,T)$ and $T\in (t_0,\tme)$,
  and in order to appropriately estimate the rightmost summand herein, we similarly employ regularity features of the
  Dirichlet heat semigroup $(e^{-t\eps A_2})_{t\ge 0}$, where $A_2$ realizes $-\Del$ in $L^p(\Om;\R^n)$ in the domain
  $D(A_2):=\{\psi\in W^{2,p}(\Om;\R^n) \ | \ \psi=0 \mbox{ on } \pO\}$, and thereby obtain $c_9=c_9(\al,\eps)>0$ such that
  \bas
	\|\ueps(\cdot,s)\|_{W^{1,p}(\Om)}
	&=& \bigg\| e^{-(t-t_0)\eps A_2} \ueps(\cdot,t_0)
	+ \int_{t_0}^t e^{-(t-s)\eps A_2} \veps(\cdot,\sig) d\sig \bigg\|_{W^{1,p}(\Om)} \\
	&\le& c_9\|\ueps(\cdot,t_0)\|_{W^{1,p}(\Om)}
	+ c_9 \int_{t_0}^t \|\veps(\cdot,\sig)\|_{W^{1,p}(\Om)} d\sig
	\qquad \mbox{for all } t\in (t_0,\tme).
  \eas
  Using that $\ueps(\cdot,t_0)\in C^2(\bom;\R^n)$, we may thus once more draw on (\ref{55.3}) to find $c_{10}=c_{10}(\al,\eps)>0$
  such that
  \bas
	\|\ueps(\cdot,s)\|_{W^{1,p}(\Om)}
	\le c_{10} + c_{10} M^{\vt_3}(T)
	\qquad \mbox{for all $s\in (t_0,T)$ and } T\in (t_0,\tme),
  \eas
  whence again relying on the inequality $\al'+\frac{1}{4}<1$, from (\ref{55.4}) we obtain $c_{11}=c_{11}(\al,\eps)>0$ such that
  \bas
	\|A^\al\veps(\cdot,t)\|_{L^p(\Om)}
	\le c_{11} + c_{11} M^{\vt_3}(T)
	\qquad \mbox{for all $t\in (t_0,T)$ and } T\in (t_0,\tme).
  \eas
  Thus,
  \bas
	M(T)
	\le c_{11} + c_{11} M^{\vt_3}(T)
	\qquad \mbox{for all } T\in (t_0,\tme),
  \eas
  which thanks to the inequality $\vt_3<1$ ensures that, indeed, $\sup_{t\in (t_0,\tme)} M(T)<\infty$.
\qed
This implies that, in fact, for any such $\eps$ the second alternative in the extensibility criterion in (\ref{ext}) cannot occur:
\begin{lem}\label{lem56}
  For each $\eps\in (0,1)$, the solution of (\ref{0eps}) from Lemma \ref{lem_loc} is global in time; that is, we have $\tme=\infty$.
\end{lem}
\proof
  If $\tme$ was finite for some $\eps\in (0,1)$, then
  fixing $\eta\in (0,\frac{1}{2})$
  we could choose an arbitrary $\al\in (\frac{1}{2},\frac{3}{4})$ and any
  $p\ge 2$ such that $4\al-\frac{n}{p}>2+2\eta$, and that thus the fractional powers in Lemma \ref{lem55} have the property that
  $D(A^\al) \hra W^{2+2\eta,\infty}(\Om;\R^n)$
(\cite[Thm.~1.6.1]{henry}).

  An application of (\ref{55.1}) would therefore yield $c_1>0$ such that, again with $t_0:=\frac{1}{2}\tm$,
  \be{56.2}
	\|\veps(\cdot,t)\|_{W^{2+2\eta,\infty}(\Om)} \le c_1
	\qquad \mbox{for all } t\in (t_0,\tme),
  \ee
  which by boundedness of $\Gamma$ on $[0,\infty)$ particularly implies that we could find $c_2>0$ such that
  $\heps:=\lan\Gamma(\Teps):\nas\veps,\nas\veps\ran$ would satisfy
  \be{56.3}
	|\heps(x,t)| \le c_2
	\qquad \mbox{for all $x\in\Om$ and } t\in (t_0,\tme).
  \ee
  Therefore, if we pick $\beta\in (\frac{1}{2},1)$ and $q>1$ large such that $2\beta-\frac{n}{q}>1+\eta$, then letting
  $B$ represent the Dirichlet Laplacian $-\Del$ in $L^q(\Om)$ and noting that $D(B^\beta) \hra W^{1+\eta,\infty}(\Om)$
  (\cite[Thm.~1.6.1]{henry}), we could employ
  standard heat semigroup estimates to obtain $c_3>0$, $c_4>0$ and $c_5>0$ fulflilling
  \bea{56.4}
	\|\Teps(\cdot,t)\|_{W^{1+\eta,\infty}(\Om)}
	&\le& c_3 \|B^\beta \Teps(\cdot,t)\|_{L^q(\Om)} \nn\\
	&=& c_3 \bigg\| B^\beta e^{-(t-t_0)B} \Teps(\cdot,t_0)
		+ \int_{t_0}^t B^\beta e^{-(t-s)B} \heps(\cdot,s) ds \bigg\|_{L^q(\Om)} \nn\\
	&\le& c_4 \|B^\beta \Teps(\cdot,t_0)\|_{L^q(\Om)}
	+ c_4 \int_{t_0}^t (t-s)^{-\beta} \|\heps(\cdot,s)\|_{L^\infty(\Om)} ds \nn\\
	&\le& c_4 \|B^\beta \Teps(\cdot,t_0)\|_{L^q(\Om)}
	+ c_2 c_4 \int_{t_0}^t (t-s)^{-\beta} ds \nn\\[2mm]
	&\le& c_5
	\qquad \mbox{for all } t\in (t_0,\tme).
 \eea
  As a combination of (\ref{56.2}) with the second equation in (\ref{0eps})
  would similarly provide $c_6>0$ satisfying
  \bas
	\|\ueps(\cdot,t)\|_{W^{1+\eta,\infty}(\Om)} \le c_6
	\qquad \mbox{for all } t\in (t_0,\tme),
  \eas
  the boundedness properties in (\ref{56.2}) and (\ref{56.4}) would contradict (\ref{ext}), however.
\qed
\mysection{Regularity properties of \tops{$\Teps$}{T eps}}
Returning to the problem of identifying $\eps$-independent properties of solutions to (\ref{0eps}), in this section we focus
on the key quantity $\Teps$ for which Lemma \ref{lem4}, also applied to $\psi\equiv 1$ here,
entails the following immediate consequence.
\begin{lem}\label{lem6}
  For each $T>0$ there exists $C(T)>0$ such that
  \be{6.1}
	\io \Teps(\cdot,t) \le C(T)
	\qquad \mbox{for all $t\in (0,T)$ and } \eps\in (0,1)
  \ee
  and
  \be{6.2}
	\int_0^T \int_{\pO} \Big| \frac{\pa\Teps}{\pa\nu}\Big| \le C(T)
	\qquad \mbox{for all }\eps\in (0,1).
  \ee
\end{lem}
\proof
  Upon applying Lemma \ref{lem4} to $\psi\equiv 1$, we see that
  \be{6.3}
	\frac{d}{dt} \io \Teps - D \int_{\pO} \frac{\pa\Teps}{\pa\nu}
	= \io \lan\Gamma(\Teps):\nas\veps,\nas\veps\ran
	\qquad \mbox{for all $t>0$ and } \eps\in (0,1),
  \ee
  where thanks to the boundedness of $\Gamma$ on $[0,\infty)$, with some $c_1>0$ we have
  \bas
	\io \lan\Gamma(\Teps):\nas\veps,\nas\veps\ran
	\le c_1 \io |\na\veps|^2
	\qquad \mbox{for all $t>0$ and } \eps\in (0,1).
  \eas
  Moreover, combining the nonnegativity of $\Teps$ with the identity $\Teps|_{\pO\times (0,\infty)}=0$ for $\eps\in (0,1)$,
  we particularly find that $\frac{\pa\Teps}{\pa\nu} \le 0$ on $\pO\times (0,\infty)$ for all $\eps\in (0,1)$, so that
  integrating (\ref{6.3}) shows that
  \bas
	\io \Teps(\cdot,t) + D \int_0^t \int_{\pO} \Big| \frac{\pa\Teps}{\pa\nu}\Big|
	\le \io \Theta_{0\eps}
	+ c_1 \int_0^t \io |\na\veps|^2
	\qquad \mbox{for all $t>0$ and } \eps\in (0,1).
  \eas
  In view of (\ref{ie}) and (\ref{5.4}), this implies both (\ref{6.1}) and (\ref{6.2}).
\qed
Based on (\ref{6.2}), we can appropriately control boundary integrals appearing in another testing procedure which is independent
from that in Lemma \ref{lem4} and, unlike the latter, capable of providing some spatially global information about temperature
gradients.
\begin{lem}\label{lem7}
  Let $p\in (0,1)$. Then for all $T>0$ there exists $C(p,T)>0$ such that
  \be{7.1}
	\int_0^T \io (\Teps+1)^{p-2} |\na\Teps|^2 \le C(p,T)
	\qquad \mbox{for all } \eps\in (0,1).
  \ee
\end{lem}
\proof
  Once more explicitly using the third equation in (\ref{0eps}), we obtain that
  \bas
	\frac{1}{p} \frac{d}{dt} \io (\Teps+1)^p
	&=& \io (\Teps+1)^{p-1} \Big\{ D\Del\Teps + \lan\Gamma(\Teps):\nas \veps,\nas \veps\ran\Big\} \\
	&\ge& D \io (\Teps+1)^{p-1} \Del\Teps
	\qquad \mbox{for all $t>0$ and } \eps\in (0,1),
  \eas
  because $\lan\Gamma(\Teps):\nas\veps,\nas\veps\ran\ge 0$ in $\Om\times (0,\infty)$ for all $\eps\in (0,1)$ due to
  (\ref{Gamma_pos}).
  Since an integration by parts shows that
  for all $t>0$ and $\eps\in (0,1)$,
  thanks to the identity
  $(\Teps+1)^{p-1}|_{\pO}=1$ we have
  \bas
	D \io (\Teps+1)^{p-1} \Del\Teps
	&=& (1-p)D\io (\Teps+1)^{p-2} |\na\Teps|^2
	+ D\int_{\pO} (\Teps+1)^{p-1} \frac{\pa\Teps}{\pa\nu} \\
	&\ge & (1-p)D\io (\Teps+1)^{p-2} |\na\Teps|^2
	- D\int_{\pO} \Big|\frac{\pa\Teps}{\pa\nu}\Big|,
  \eas
  this entails that
  \bea{7.2}
	(1-p)D \int_0^T \io (\Teps+1)^{p-2} |\na\Teps|^2
	&\le& \frac{1}{p} \io \Big(\Teps(\cdot,T)+1\Big)^p
	- \frac{1}{p} \io (\Theta_{0\eps}+1)^p
	+ D\int_0^T \int_{\pO} \Big| \frac{\pa\Teps}{\pa\nu}\Big| \nn\\
	&\le& \frac{1}{p} \io \Big(\Teps(\cdot,T)+1\Big)^p
	+ D\int_0^T \int_{\pO} \Big| \frac{\pa\Teps}{\pa\nu}\Big|
  \eea
  for all $T>0$ and $\eps\in (0,1)$.
  Using Young's inequality in estimating
  \bas
	\frac{1}{p} \io \Big(\Teps(\cdot,T)+1\Big)^p
	\le \frac{1}{p} \io \Big(\Teps(\cdot,T)+1\Big) + \frac{|\Om|}{p}
	\qquad \mbox{for all $T>0$ and } \eps\in (0,1),
  \eas
  in light of (\ref{6.1}) and (\ref{6.2}) we infer (\ref{7.1}) from (\ref{7.2}).
\qed
Through straightforward Gagliardo-Nirenberg interpolation, this implies bounds for $\Teps$ in some
reflexive Lebesgue spaces.
\begin{lem}\label{lem8}
  Let $q\in (1,\frac{n+2}{n})$.
  Then for each $T>0$ there exists $C(q,T)>0$ such that
  \be{8.1}
	\int_0^T \io (\Teps+1)^q \le C(q,T)
	\qquad \mbox{for all } \eps\in (0,1).
  \ee
\end{lem}
\proof
  By boundedness of $\Om$, we may assume without loss of generality that $q>\frac{2}{n}$, and that thus our assumption implies
  that $p\equiv p(q):=q-\frac{2}{n}$ satisfies $p\in (0,1)$.
  An application of the Gagliardo-Nirenberg inequality then yields $c_1=c_1(q)>0$ such that
  \be{8.2}
	\io |\psi|^\frac{2q}{p}
	\le c_1\|\na\psi\|_{L^2(\Om)}^{\frac{2q}{p}\vt} \|\psi\|_{L^\frac{2}{p}(\Om)}^{\frac{2q}{p}(1-\vt)}
	+ c_1\|\psi\|_{L^\frac{2}{p}(\Om)}^\frac{2q}{p}
	\qquad \mbox{for all } \psi\in W^{1,2}(\Om;\R),
  \ee
  where by definition of $p$, the number $\vt\equiv \vt(q):=(\frac{np}{2}-\frac{np}{2q})/(1-\frac{n}{2}+\frac{np}{2})$ satisfies
  \bas
	\frac{2q}{p}\vt
	= \frac{nq-n}{1-\frac{n}{2} + \frac{np}{2}}
	= \frac{nq-n}{1-\frac{n}{2}+\frac{n}{2}(q-\frac{2}{n})}
	= \frac{nq-n}{-\frac{n}{2}+\frac{nq}{2}}=2.
  \eas
  Therefore, (\ref{8.2}) implies that
  \bea{8.3}
	\io (\Teps+1)^q
	&=& \io \Big\{ (\Teps+1)^\frac{p}{2}\Big\}^\frac{2q}{p} \nn\\
	&\le& c_1 \bigg\{ \io \big|\na (\Teps+1)^\frac{p}{2}\big|^2 \bigg\} \bigg\{ \io (\Teps+1)\bigg\}^\frac{2}{n}
	+ c_1\bigg\{ \io (\Teps+1)\bigg\}^q \nn\\
	&\le& c_1 c_2^\frac{2}{n} \io \big|\na (\Teps+1)^\frac{p}{2}\big|^2
	+ c_1 c_2^q
	\qquad \mbox{for all $t\in (0,T)$ and } \eps\in (0,1),
  \eea
  where $c_2\equiv c_2(T):=\sup_{\eps\in (0,1)} \sup_{t\in (0,T)} \io (\Teps(\cdot,t)+1)$ is finite due to Lemma \ref{lem6}.
  Since $|\na (\Teps+1)^\frac{p}{2}|^2 = \frac{p^2}{4} (\Teps+1)^{p-2} |\na\Teps|^2$ for all $\eps\in (0,1)$, integrating
  (\ref{8.3}) over $t\in (0,T)$ and using Lemma \ref{lem7} we arrive at (\ref{8.1}).
\qed
One further standard interpolation step yields bounds for $\na\Teps$ which in contrast to those in Lemma \ref{lem7}
do no longer contain weight functions.
\begin{lem}\label{lem9}
  If $r\in (1,\frac{n+2}{n+1})$,
  then given any $T>0$ one can find $C(r,T)>0$ such that
  \be{9.1}
	\int_0^T \io |\na\Teps|^r \le C(r,T)
	\qquad \mbox{for all } \eps\in (0,1).
  \ee
\end{lem}
\proof
  Since $r<\frac{n+2}{n+1}$ and hence $(3n+2)r-2(n+2)<nr$, we can fix $p=p(r)\in (0,1)$ suitably close to $1$ such that
  \bas
	p>\frac{(3n+2)r-2(n+2)}{nr}.
  \eas
  As thus $2-p<\frac{n+2}{n}\cdot \frac{2-r}{r}$, letting
  $q\equiv q(r):=\frac{(2-p)r}{2-r}$ defines a number $q\in (1,\frac{n+2}{n})$, and since Young's inequality says that
  according to this definition of $q$ we have
  \bas
	\int_0^T \io |\na\Teps|^r
	&=& \int_0^T \io \Big\{ (\Teps+1)^{p-2} |\na\Teps|^2 \Big\}^\frac{r}{2} (\Teps+1)^\frac{(2-p)r}{2} \\
	&\le& \int_0^T \io (\Teps+1)^{p-2} |\na\Teps|^2
	+ \int_0^T \io (\Teps+1)^q
	\qquad \mbox{for all $T>0$ and } \eps\in (0,1),
  \eas
  the claim is a consequence of Lemma \ref{lem7} and Lemma \ref{lem8}.
\qed
In a natural manner, Lemma \ref{lem8} and Lemma \ref{lem5} imply some information on regularity of $\Theta_{\eps t}$
that will be used to extract pointwise a.e.~convergent subsequences by means of an Aubin-Lions lemma (Lemma \ref{lem11}).
\begin{lem}\label{lem10}
  Let $s>\frac{n+2}{2}$. Then for all $T>0$ there exists $C(s,T)>0$ such that
  \be{10.1}
	\int_0^T \big\|\Theta_{\eps t}(\cdot,t)\big\|_{(W_0^{2,s}(\Om;\R))^\star} dt \le C(s,T)
	\qquad \mbox{for all } \eps\in (0,1).
  \ee
\end{lem}
\proof
  For fixed $\psi\in C_0^\infty(\Om;\R)$, again due to the boundedness of $\Gamma$ on $[0,\infty)$ the identity in (\ref{4.1})
  shows that with some $c_1>0$,
  \bas
	\bigg| \io \Theta_{\eps t} \psi \bigg|
	&=& \bigg| D\io \Teps \Del\psi
	+ \io \lan\Gamma(\Teps):\nas\veps,\nas\veps\ran \psi \bigg| \\
	&\le& D \io \Teps |\Del\psi|
	+ c_1 \io |\na\veps|^2 |\psi| \\
	&\le& D\|\Teps\|_{L^\frac{s}{s-1}(\Om)} \|\Del\psi\|_{L^s(\Om)}
	+ c_1\|\na\veps\|_{L^2(\Om)}^2 \|\psi\|_{L^\infty(\Om)}
	\qquad \mbox{for all $t>0$ and } \eps\in (0,1).
  \eas
  Since the inequality $s>\frac{n+2}{2}>\frac{n}{2}$ particularly ensures that $W^{2,s}(\Om;\R)$ is continuously embedded into
  $L^\infty(\Om;\R)$, by definition of the norm in $(W_0^{2,s}(\Om;\R))^\star$ this entails the existence of $c_2=c_2(s)>0$
  fulfilling
  \bas
	\|\Theta_{\eps t}\|_{(W_0^{2,s}(\Om;\R))^\star} \le D \|\Teps\|_{L^\frac{s}{s-1}(\Om)}
	+ c_1c_2\|\na\veps\|_{L^2(\Om)}^2
	\qquad \mbox{for all $t>0$ and } \eps\in (0,1).
  \eas
  Young's inequality thus implies that
  for all $T>0$ and $\eps\in (0,1)$,
  \bas
	\int_0^T \big\| \Theta_{\eps t}(\cdot,t)\big\|_{(W_0^{2,s}(\Om;\R))^\star} dt
	&\le& D \int_0^T \|\Teps(\cdot,t)\|_{L^\frac{s}{s-1}(\Om)} dt
	+ c_1c_2 \int_0^T \|\na\veps(\cdot,t)\|_{L^2(\Om)}^2 dt \\
	&\le& D \int_0^T \io \Teps^\frac{s}{s-1}
	+ D T
	+ c_1c_2 \int_0^T \io |\na\veps|^2,
  \eas
  so that the claim results upon employing (\ref{5.1}) and Lemma \ref{lem8}, and observing that
  $\frac{s}{s-1}=(1-\frac{1}{s})^{-1} < (1-\frac{2}{n+2})^{-1} = \frac{n+2}{n}$ thanks to our hypothesis.
\qed
\mysection{Exploiting weak lower semicontinuity of \tops{$L^2$}{L²} norms}\label{sect_lsc}
This section collects some consequences of lower semicontinuity of norms in $L^2$ spaces with respect to weak convergence,
arranged here in such a way that expressions of the form in (\ref{quad}) can adequately be coped with in several different
particular situations arising below.\abs
Let us first record an essentially well-known consequence of Lebesgue's theorem for products of a.e.~convergent and $L^2$-weakly
convergent sequences.
This will be used not only in the subsequent Lemma \ref{lem1}, but later on also in our extraction of
convergent subsequences of $((\veps,\ueps,\Teps))_{\eps\in (0,1)}$ (Lemma \ref{lem11}) and in the derivation of (\ref{wF})
(Lemma \ref{lem12}).
\begin{lem}\label{lem19}
  Let $N\ge 1$ and $G\subset \R^N$ be measurable with $|G|<\infty$,
  and suppose that $(f_j)_{j\in\N} \subset L^\infty(G)$,
  $(g_j)_{j\in\N} \subset L^2(G)$, $f\in L^\infty(G)$ and $g\in L^2(G)$ are such that
  \be{19.1}
	\sup_{j\in\N} \|f_j\|_{L^\infty(G)} < \infty,
  \ee
  and that as $j\to\infty$ we have
  \be{19.2}
	f_j \to f
	\qquad \mbox{a.e.~in } G
  \ee
  and
  \be{19.3}
	g_j \wto g
	\qquad \mbox{in } L^2(G).
  \ee
  Then
  \be{19.4}
	f_j g_j \wto fg
	\quad \mbox{in } L^2(G)
	\qquad \mbox{as } j\to\infty.
  \ee
\end{lem}
\proof
  If the claim was false, then since (\ref{19.1}) and (\ref{19.3}) particularly assert that $(f_j g_j)_{j\in\N}$ is bounded
  in $L^2(G)$, we could find a subsequence $(j_k)_{k\in\N}$ along which $f_{j_k} g_{j_k} \wto h$ would hold in $L^2(G)$
  as $k\to\infty$
  with some $h\in L^2(G)$ satisfying $|\{h\ne fg\}| >0$. As $|G|$ is finite, this would especially imply that
  $f_{j_k} g_{j_k} \wto h$ in $L^1(G)$ as $k\to\infty$.
  But uniform boundedness of $(f_j)_{j\in ℕ}$ according to \eqref{19.1} together with \eqref{19.2} shows that $f_j\to f$ in $L^2(G)$ by Lebesgue's dominated convergence theorem. Accordingly, $f_jg_j\wto g$ in $L^1(G)$ due to \eqref{19.3}, and hence we could infer that $h=fg$, which is absurd.
\qed
\begin{lem}\label{lem:root}
   Let $\Beta=(\Beta_{ijkl})_{i,j,k,l\in\{1,...,n\}} \in C^0([0,\infty);\R^{n\times n\times n \times n}) \cap
  L^\infty([0,\infty);\R^{n\times n\times n \times n})$ be such that
  \be{1.1}
	\Beta_{ijkl}(\xi)=\Beta_{klij}(\xi)
	\qquad \mbox{for all $\xi\ge 0$ and } (i,j,k,l)\in\{1,...,n\}^4,
  \ee
  and that there exists $K_\Beta>0$ such that
  \be{1.2}
	\lan \Beta(\xi):X,X\ran \ge K_\Beta |X|^2
	\qquad \mbox{for all $\xi\ge 0$ and } X\in\R^{n\times n}.
  \ee
 Then there exists $\sqrt{\Beta}\in C^0([0,\infty);\R^{n\times n\times n \times n}) \cap
  L^\infty([0,\infty);\R^{n\times n\times n \times n})$
  such that
  \be{1.77}
	\lan\Beta(\xi):W,W\ran = \lan \sqrt{\Beta(\xi)}:W,\sqrt{\Beta(\xi)}:W\ran
	\qquad \mbox{for all $\xi\ge 0$ and } W\in \R^{n\times n}.
  \ee
\end{lem}
\proof
 For every $ξ\in[0,∞)$ interpreting $\Beta(ξ)$ as (by \eqref{1.1}) symmetric and (by \eqref{1.2}) positive definite $n^2\times n^2$ matrix, we take $\sqrt{\Beta}(ξ)$ to be its unique positive definite square root \cite[Cor. 1.30]{higham}, which in particular satisfies \eqref{1.77}. Then $ξ\mapsto \sqrt{\Beta}(ξ)$ is continuous by continuity of $\Beta$ and \cite[Thm. 6.12]{higham}. Boundedness of $\sqrt{\Beta}$ follows from $\Beta\in L^\infty([0,\infty);\R^{(n\times n)\times (n \times n)})$ and is most easily seen in the spectral norm, as $\norm[σ]{\sqrt{B}(ξ)}=\sqrt{\norm[σ]{B(ξ)}}$.
\qed

For bounded and uniformly positive definite symmetric tensor-valued functions, we can use Lemma \ref{lem19} together with the
lower semicontinuity property under consideration to obtain the following basic property.
\begin{lem}\label{lem1}
  Let $\Beta=(\Beta_{ijkl})_{i,j,k,l\in\{1,...,n\}} \in C^0([0,\infty);\R^{n\times n\times n \times n}) \cap
  L^\infty([0,\infty);\R^{n\times n\times n \times n})$ satisfy \eqref{1.1} and let there be $K_\Beta>0$ such that \eqref{1.2} holds.
  Then whenever $T>0$ as well as $\vp\in L^\infty(\Om\times (0,T);\R)$,
  $(W_j)_{j\in\N} \subset L^2(\Om\times (0,T);\R^{n\times n})$ and
  the measurable functions $z_j:\Om\times (0,T)\to\R$, $j\in\N$, are such that
  \be{1.3}
	\vp\ge 0
	\qquad \mbox{a.e.~in } \Om\times (0,T)
  \ee
  and
  \be{1.4}
	z_j\ge 0
	\quad \mbox{a.e.~in } \Om\times (0,T)
	\qquad \mbox{for all } j\in\N,
  \ee
  and that as $j\to\infty$ we have
  \be{1.5}
	W_j \wto W
	\qquad \mbox{in } L^2(\Om\times (0,T);\R^{n\times n})
  \ee
  and
  \be{1.6}
	z_j \to z
	\qquad \mbox{a.e.~in } \Om\times (0,T)
  \ee
  with some $W\in L^2(\Om\times (0,T);\R^{n\times n})$ and some measurable $z:\Om\times (0,T)\to \R$, it follows that
  \be{1.7}
	\int_0^T \io \lan \Beta(z):W,W\ran \vp
	\le \liminf_{j\to\infty} \int_0^T \io \lan \Beta(z_j):W_j,W_j\ran \vp.
  \ee
\end{lem}
\proof
  According to (\ref{1.1}), (\ref{1.2}) and Lemma~\ref{lem:root},
  there exists $\sqrt{\Beta}\in C^0([0,\infty);\R^{n\times n\times n \times n}) \cap
  L^\infty([0,\infty);\R^{n\times n\times n \times n})$
  such that \eqref{1.77} holds.

  Thus, if for $j\in\N$ we let
  \be{1.78}
	\rho_j:=\sqrt{\vp} \sqrt{\Beta(z_j)} : W_j,
  \ee
  then we obtain a sequence $(\rho_j)_{j\in\N} \subset L^2(\Om\times (0,T);\R^{n\times n})$ which due to Lemma \ref{lem19}
  and (\ref{1.5}) satisfies
  \be{1.87}
	\rho_j \wto \rho:=\sqrt{\vp} \sqrt{\Beta(z)}:W
	\quad \mbox{in } L^2(\Om\times (0,T);\R^{n\times n})
	\qquad \mbox{as } j\to\infty,
  \ee
  because $\sup_{j\in\N} \|\sqrt{\Beta(z_j)}\|_{L^\infty(\Om\times (0,T))} < \infty$ and $\sqrt{\Beta(z_j)} \to \sqrt{\Beta(z)}$
  a.e.~in $\Om\times (0,T)$ by boundedness and continuity of $\sqrt{\Beta}$, and by (\ref{1.6}).
  Thanks to the lower semicontinuity of the norm in $L^2(\Om\times (0,T);\R^{n\times n})$ with respect to weak convergence,
  from (\ref{1.87}) we infer that, in line with (\ref{1.77}),
  \bas
	\int_0^T \io \lan \Beta(z):W,W\ran \vp
	= \int_0^T \io |\rho|^2
	\le \liminf_{j\to\infty} \int_0^T \io |\rho_j|^2
	= \liminf_{j\to\infty} \int_0^T \io \lan \Beta(z_j):W_j,W_j\ran \vp,
  \eas
  as intended.
\qed

Of crucial importance for our reasoning is now the following consequence of Lemma~\ref{lem1}.
\begin{lem}\label{lem102}
  Let $\beta=(\beta_{ijkl})_{i,j,k,l\in\{1,...,n\}} \in C^0([0,\infty);\R^{n\times n\times n \times n}) \cap
  L^\infty([0,\infty);\R^{n\times n\times n \times n})$ be such that
  \be{102.01}
	\beta_{ijkl}(\xi)=\beta_{klij}(\xi)
	\qquad \mbox{for all $\xi\ge 0$ and } (i,j,k,l)\in\{1,...,n\}^4,
  \ee
  Then there exists $\eta_1=\eta_1(\beta)>0$ with the property that if $\eta\in (0,\eta_1)$ and $T>0$,
  if $\vp\in L^\infty(\Om\times (0,T);\R)$ is nonnegative, and if
  $(w_j)_{j\in\N} \subset L^2((0,T);W_0^{1,2}(\Om;\R^n))$,
  $(z_j)_{j\in\N} \subset L^1(\Om\times (0,T);\R)$,
  $w\in L^2((0,T);W_0^{1,2}(\Om;\R^n))$,
  and $z \in L^1(\Om\times (0,T);\R)$
  are such that
  \be{102.2}
	z_j\ge 0
	\mbox{ a.e.~in } \Om\times (0,T)
	\quad \mbox{for all } j\in\N
	\qquad \mbox{and} \qquad
	z_j \to z
	\mbox{ a.e.~in } \Om\times (0,T)
	\qquad \mbox{as } j\to\infty,
  \ee
  and that
  \be{102.1}
	w_j \wto w
	\quad \mbox{in } L^2((0,T);W_0^{1,2}(\Om;\R^n))
	\qquad \mbox{as } j\to\infty,
  \ee
  then
  \bea{102.3}
	& & \hs{-30mm}
	\int_0^T \io |\na w|^2 \vp
	- \eta \int_0^T \io \lan \beta(z):\nas w,\nas w\ran \vp \nn\\
	&\le& \liminf_{j\to\infty} \bigg\{
	\int_0^T \io |\na w_j|^2 \vp
	- \eta \int_0^T \io \lan \beta(z_j):\nas w_j, \nas w_j\ran \vp \bigg\}.
  \eea
\end{lem}
\proof
  We let $(\del_{ij})_{i,j\in\{1,...,n\}}$ denote the Kronecker delta on $\{1,...,n\}^2$, and define
  \be{102.44}
	\wh{\beta}_{ijkl}:=\frac{1}{2}(\del_{ik}\del_{jl}+\del_{il}\del_{jk}),
	\qquad (i,j,k,l)\in\{1,...,n\}^4,
  \ee
  observing that then
  \bas
	\big(\wh{\beta}:X\big)_{ij}
	= \frac{1}{2} \sum_{k,l=1}^n \del_{ik} \del_{jl} X_{kl}
	+ \frac{1}{2} \sum_{k,l=1}^n \del_{il} \del_{jk} X_{kl}
	= \frac{1}{2} X_{ij} + \frac{1}{2} X_{ji}
	\qquad \mbox{for all } (i,j)\in\{1,...,n\}^2
  \eas
  and hence
  \bas
	\wh{\beta}:X = \frac{1}{2} \big(X+X^t\big)
  \eas
  for each $X=(X_{ij})_{i,j\in\{1,...,n\}} \in \R^{n\times n}$.
  As clearly
  \be{102.4}
	\wh{\beta}_{ijkl} = \wh{\beta}_{klij}
	\qquad \mbox{for all } (i,j,k,l)\in\{1,...,n\}^4,
  \ee
  by straightforward linear algebra we thus obtain that
  \bea{102.5}
	\lan \beta(\xi):\Big\{ \frac{1}{2}\big( X+X^t\big) \Big\}, \frac{1}{2}\big(X+X^t\big) \ran
	&=& \lan \beta(\xi):\big(\wh{\beta}:X\big),\wh{\beta}:X\ran \nn\\
	&=& \lan \wh{\beta}:\big\{ \beta(\xi):\big(\wh{\beta}:X\big)\big\},X\ran \nn\\
	&=& \lan \wh{\Beta}(\xi):X,X\ran
	\qquad \mbox{for all $\xi\ge 0$ and } X\in\R^{n\times n}
  \eea
  with
  \be{102.6}
	\wh{\Beta}(\xi):= \wh{\beta} \omult \big( \beta(\xi)\omult \wh{\beta}\big),
	\qquad \xi\ge 0,
  \ee
  where for $\beta^{(\iota)}=(\beta^{(\iota)}_{ijkl})_{i,j,k,l\in\{1,...,n\}}\in \R^{n\times n\times n\times n}$,
  $\iota\in \{1,2\}$, we have set
  \bas
	\big(\beta^{(1)} \omult \beta^{(2)} \big)_{ijkl}
	:= \sum_{m,m'=1}^n \beta^{(1)}_{ijmm'} \beta^{(2)}_{mm' kl},
	\qquad (i,j,k,l)\in\{1,...,n\}^4.
  \eas
  Now due to the assumed boundeness of $\beta$ and (\ref{102.44}), from (\ref{102.6}) we infer that with some
  $c_1=c_1(\beta)>0$ we have
  \be{102.7}
	\lan \wh{\Beta}(\xi):X,X\ran
	\le c_1 |X|^2
	\qquad \mbox{for all $\xi\ge 0$ and } X\in\R^{n\times n},
  \ee
  and fixing any $\eta_1=\eta_1(\beta)\in (0,\frac{1}{c_1})$ and assuming that $\eta\in (0,\eta_1)$, we thereupon let
  \be{102.8}
	\Beta(\xi):=I-\eta \wh{\Beta}(\xi),
	\qquad \xi\ge 0,
  \ee
  where $I:=(\del_{ik}\del_{jl})_{i,j,k,l\in\{1,...,n\}}$ represents the identity mapping on $\R^{n\times n}$.
  Again thanks to (\ref{102.4}), the hypothesis in (\ref{102.01}) ensures that $\Beta$ has the symmetry property in (\ref{1.1}),
  while (\ref{102.8}) along with (\ref{102.7}) guarantees that
  \bas
	\lan \Beta(\xi):X,X\ran
	= |X|^2 - \eta \lan \wh{\Beta}(\xi):X,X\ran
	\ge c_2 |X|^2
	\qquad \mbox{for all $\xi\ge 0$ and } X\in\R^{n\times n},
  \eas
  with $c_2\equiv c_2(\beta):=1-\eta_1 c_1$ being positive due to our selection of $\eta_1$.
  The lemma thereby becomes a consequence of \eqref{102.1} and Lemma \ref{lem1}, according to which, namely, in line with (\ref{102.8})
  and (\ref{102.6}) it follows that
  \bas
	& & \hs{-20mm}
	\int_0^T \io |\na w|^2 \vp
	- \eta \int_0^T \io \lan \beta(z):\nas w,\nas w\ran\vp  \\
	&=& \int_0^T \io \lan \Beta(z):\na w,\na w\ran \vp \\
	&\le& \liminf_{j\to\infty} \int_0^T \io \lan \Beta(z_j):\na w_j,\na w_j\ran \vp \\
	&=& \liminf_{j\to\infty} \bigg\{
	\int_0^T \io |\na w_j|^2 \vp
	- \eta \int_0^T \io \lan \beta(z_j):\nas w_j,\nas w_j\ran \vp \bigg\},
  \eas
  as claimed.
\qed
In preparation for a relative of Lemma \ref{lem102} involving differences with ordering opposite to those in (\ref{102.3}),
let us state the following 
observation which essentially relies on vanishing boundary values.
\begin{lem}\label{lem1033}
  Let $w\in W_0^{1,2}(\Om;\R^n)$ and $\psi\in C^1(\bom;\R)$.
  Then
  \be{1033.1}
	\io |\nas w|^2 \psi
	= \frac{1}{2} \io |\na w|^2 \psi
	+ \frac{1}{2} \io |\div w|^2 \psi
	+ \frac{1}{2} \io (\div w) (w\cdot\na\psi)
	- \frac{1}{2} \io w\cdot (\na w\cdot\na\psi).
  \ee
\end{lem}
\proof
  For $w=(w_1,...,w_n)\in C_0^\infty(\Om;\R^n)$ and $\psi\in C^1(\bom;\R)$, in the identity
  \bas
	\io |\nas w|^2 \psi
	&=& \frac{1}{4} \sum_{i,j=1}^n \io (\pa_j w_i + \pa_i w_j)^2 \psi \\
	&=& \frac{1}{4} \bigg\{ \sum_{i,j=1}^n \io (\pa_j w_i)^2 \psi
	+ \sum_{i,j=1}^n \io (\pa_i w_j)^2 \psi
	+ 2\sum_{i,j=1}^n \io \pa_j w_i \pa_i w_j \psi
	\bigg\} \\
	&=& \frac{1}{2} \sum_{i,j=1}^n\io (\pa_j w_i)^2 \psi
	+ \frac{1}{2} \sum_{i,j=1}^n \io\pa_i w_j \pa_j w_i \psi \\
	&=& \frac{1}{2} \io |\na w|^2 \psi
	+ \frac{1}{2} \sum_{i,j=1}^n\io \pa_i w_j \pa_j w_i \psi
  \eas
  we may twice integrate by parts to see that
  \bas
	\frac{1}{2} \sum_{i,j=1}^n\io \pa_i w_j \pa_j w_i \psi
	&=& - \frac{1}{2} \sum_{i,j=1}^n\io w_j \pa_{ij} w_i \psi
	- \frac{1}{2} \sum_{i,j=1}^n\io w_j \pa_j w_i \pa_i\psi \\
	&=& \frac{1}{2} \sum_{i,j=1}^n\io \pa_j w_j \pa_i w_i \psi
	+ \frac{1}{2} \sum_{i,j=1}^n\io w_j \pa_i w_i \pa_j \psi
	- \frac{1}{2} \sum_{i,j=1}^n\io w_j \pa_j w_i \pa_i\psi \\
	&=& \frac{1}{2} \io |\div w|^2 \psi
	+ \frac{1}{2} \io (\div w) (w\cdot\na\psi)
	- \frac{1}{2} \io w\cdot (\na w\cdot\na\psi),
  \eas
  so that (\ref{1033.1}) immediately follows by means of a completion argument.
\qed
In fact, making appopriate use of this we can complement Lemma \ref{lem102} by means of a different argument as follows.
\begin{lem}\label{lem103}
  Suppose that $\beta=(\beta_{ijkl})_{i,j,k,l\in\{1,...,n\}} \in C^0([0,\infty);\R^{n\times n\times n \times n}) \cap
  L^\infty([0,\infty);\R^{n\times n\times n \times n})$ satisfies \eqref{102.01} and
  \be{103.1}
	\lan \beta(\xi):X,X\ran \ge K_\beta |X|^2
	\qquad \mbox{for all $\xi\ge 0$ and } X\in\R^{n\times n}
  \ee
  with some $K_\beta>0$.
  Then there exists $\eta_2=\eta_2(\beta)>0$ such that assuming that $T>0$, that
  $0\le \vp\in L^\infty(\Om\times (0,T);\R)$, that
  $(w_j)_{j\in\N} \subset L^2((0,T);W_0^{1,2}(\Om;\R^n))$,
  $(z_j)_{j\in\N} \subset L^1(\Om\times (0,T);\R)$,
  $w\in L^2((0,T);W_0^{1,2}(\Om;\R^n))$,
  and $z \in L^1(\Om\times (0,T);\R)$ satisfy (\ref{102.2}), (\ref{102.1}) as well as
  \be{103.2}
	w_j \to w
	\quad \mbox{in } L^2(\Om\times (0,T);\R^n)
	\qquad \mbox{as } j\to\infty,
  \ee
  one can conclude that
  \bea{103.3}
	& & \hs{-40mm}
	\int_0^T \io \lan \beta(z):\nas w,\nas w\ran \vp
	- \eta \int_0^T \io |\na w|^2 \vp \nn\\
	&\le& \liminf_{j\to\infty} \bigg\{
	\int_0^T \io \lan \beta(z_j):\nas w_j, \nas w_j\ran \vp
	- \eta 	\int_0^T \io |\na w_j|^2 \vp \bigg\}
  \eea
  for all $\eta\in (0,\eta_2)$.
\end{lem}
\proof
  We fix any $\eta_2=\eta_2(\beta)\in (0,\frac{1}{2}K_\beta)$, and for $\eta\in (0,\eta_2)$ we let
  \bas
	\Beta(\xi):=\beta(\xi)-2\eta I,
	\qquad \xi\ge 0,
  \eas
  where again $I:=(\del_{ik}\del_{jl})_{i,j,k,l\in\{1,...,n\}}\in\R^{n\times n\times n \times n}$.
  Then $\Beta$ satisfies (\ref{1.1}) due to (\ref{102.01}), whereas the requirement in (\ref{103.1}) guarantees that
  writing $c_1\equiv c_1(\beta):=K_\beta-2\eta_2>0$ we have
  \bas
	\lan \Beta(\xi):X,X\ran
	= \lan \beta(\xi):X,X\ran
	- 2\eta |X|^2
	\ge c_1 |X|^2
	\qquad \mbox{for all $\xi\ge 0$ and } X\in\R^{n\times n}.
  \eas
  Since $\nas w_j\wto \nas w$ in $L^2(\Om\times (0,T);\R^{n\times n})$ as $j\to\infty$, we may therefore employ Lemma \ref{lem1}
  to infer that
  \bea{103.4}
	& & \hs{-30mm}
	\int_0^T \io \lan \beta(z):\nas w,\nas w\ran \vp
	- 2\eta \int_0^T \io |\nas w|^2 \vp \nn\\
	&=& \int_0^T \io \lan \Beta(z):\nas w,\nas w\ran \vp \nn\\
	&\le& \liminf_{j\to\infty} \int_0^T \io \lan \Beta(z_j):\nas w_j,\nas w_j\ran \vp \nn\\
	&=& \liminf_{j\to\infty} \bigg\{
	\int_0^T \io \lan \beta(z_j):\nas w_j,\nas w_j\ran \vp
	- 2\eta \int_0^T \io |\nas w_j|^2 \vp
	\bigg\}.
  \eea
  To appropriately complement this, we now utilize Lemma \ref{lem1033} to confirm that
  for all $j\in\N$,
  \bea{103.5}
	& & \hs{-24mm}
	2\eta \int_0^T \io |\nas w_j|^2 \vp
	- \eta \int_0^T \io |\na w_j|^2 \vp \nn\\
	&=& \eta \int_0^T \io |\div w_j|^2 \vp
	+ \eta \int_0^T \io (\div w_j) (w_j\cdot\na\vp)
	- \eta \int_0^T \io w_j\cdot (\na w_j\cdot\na\vp),
  \eea
  where as a consequence of (\ref{102.1}) and the strong convergence feature in (\ref{103.2}),
  \be{103.6}
	\eta \int_0^T \io (\div w_j) (w_j\cdot\na\vp)
	\to \eta \int_0^T \io (\div w) (w\cdot\na\vp)
  \ee
  and
  \be{103.7}
	- \eta \int_0^T \io w_j\cdot (\na w_j\cdot\vp)
	\to - \eta \int_0^T \io w\cdot (\na w\cdot\na\vp)
  \ee
  as $j\to\infty$.
  Since by lower semicontinuity of the norm in $L^2(\Om\times (0,T);\R)$ we furthermore readily obtain from (\ref{102.1}) that
  \bas
	\eta \int_0^T \io |\div w|^2 \vp
	\le \liminf_{j\to\infty} \bigg\{ \eta \int_0^T \io |\div w_j|^2 \vp \bigg\},
  \eas
  once more relying on Lemma \ref{lem1033} we may combine (\ref{103.5}) with (\ref{103.6}) and (\ref{103.7}) to see that
  \bas
	& & \hs{-20mm}
	2\eta \int_0^T \io |\nas w|^2 \vp
	- \eta \int_0^T \io |\na w|^2 \vp \\
	&=& \eta \int_0^T \io |\div w|^2 \vp
	+ \eta \int_0^T \io (\div w)(w\cdot\na\vp)
	- \eta \int_0^T \io w\cdot (\na w\cdot\na\vp) \\
	&\le& \liminf_{j\to\infty} \bigg\{
	\eta \int_0^T \io |\div w_j|^2 \vp
	+ \eta \int_0^T \io (\div w_j) (w_j\cdot\na\vp)
	- \eta \int_0^T \io w_j\cdot (\na w_j\cdot\na\vp)
	\bigg\} \\
	&=& \liminf_{j\to\infty} \bigg\{
	2\eta \int_0^T \io |\nas w_j|^2 \vp
	- \eta \int_0^T \io |\na w_j|^2 \vp
	\bigg\}.
  \eas
  In conjunction with (\ref{103.4}), this establishes (\ref{103.3}) due to the basic fact that
  $\liminf_{j\to\infty} (\om_j+\wh{\om}_j) \ge \liminf_{j\to\infty} \om_j + \liminf_{j\to\infty} \wh{\om}_j$
  for bounded sequences $(\om_j)_{j\in\N}\subset\R$ and $(\wh{\om}_j)_{j\in\N}\subset\R$.
\qed
\mysection{Passing to the limit. Proof of Theorem \ref{theo13}}
Final preliminaries for our limit passage in (\ref{0eps}) derive the following information on regularity of temporal derivatives 
from Lemma \ref{lem5}.
\begin{lem}\label{lem110}
  Let $T>0$. Then there exists $C(T)>0$ such that
  \be{110.1}
	\int_0^T \big\| v_{\eps t}(\cdot,t)\big\|_{(W_0^{2,2}(\Om;\R^n))^\star}^2 dt \le C(T)
	\qquad \mbox{for all } \eps\in (0,1).
  \ee
\end{lem}
\proof
  Given $\psi\in C_0^\infty(\Om;\R^n)$ such that $\|\vp\|_{W_0^{2,2}(\Om)} \le 1$, we integrate by parts in the first equation
  from (\ref{0eps}) to see that due to the Cauchy-Schwarz inequality,
  \bas
	\bigg| \io v_{\eps t} \cdot\psi \bigg|
	&=& \bigg| - \eps \io \Del^2 \veps\cdot\psi
	+ \io \div \big(\gamma(\Teps):\nas (\veps+a\ueps)\big)\cdot\psi \bigg| \\
	&=& \bigg| - \eps \io \Del\veps\cdot\Del\psi
	- \io \lan \gamma(\Teps):\nas (\veps+a\ueps),\na\psi \ran \bigg| \\
	&\le& \eps^\frac{1}{2} \bigg\{ \eps \io |\Del\veps|^2 \bigg\}^\frac{1}{2} \|\Del\psi\|_{L^2(\Om)}
	+ \bigg\{ \io \big|\gamma(\Teps):\nas (\veps+a\ueps)\big|^2 \bigg\}^\frac{1}{2} \|\na\psi\|_{L^2(\Om)}
  \eas
  for all $t>0$ and $\eps\in (0,1)$.
  According to the boundedness of $\gamma$ and the definition of the norm in $(W_0^{2,2}(\Om;\R^n))^\star$, this implies the
  existence of $c_1>0$ such that
  \bas
	\big\| v_{\eps t}(\cdot,t)\big\|_{(W_0^{2,2}(\Om;\R^n))^\star}^2
	\le c_1 \eps \io |\Del\veps|^2
	+ c_1 \io |\na\veps|^2
	+ c_1\io |\na\ueps|^2
	\qquad \mbox{for all $t>0$ and } \eps\in (0,1),
  \eas
  which in view of (\ref{5.5}), (\ref{5.4}) and (\ref{5.2}) entails \eqref{110.1} upon an integration in time.
\qed

A similar statement holds for the derivative $u_{\eps t}$.
\begin{lem}\label{lem:ut}
  Let $T>0$. Then there exists $C(T)>0$ such that
\[
 \norm[L^\infty((0,T);(W_0^{2,2}(\Omega;ℝ^n))^\star)]{u_{\eps t}}\le C(T) \qquad \mbox{for all } \eps\in (0,1).
\]
\end{lem}
\proof
Given $ψ\in C_0^\infty(\Omega;ℝ^n)$ such that $\norm[W_0^{2,2}(\Om)]{ψ}\le 1$, we conclude from the second equation of \eqref{0eps}, integration by parts and Hölder's inequality that
\[
 \left| \io u_{\eps t}(\cdot,t) ψ \right| \le ε\io |\ueps(\cdot,t)| |Δψ| + \io \veps |ψ| \le ε\norm[\Lom2]{\ueps(\cdot,t)}+\norm[\Lom2]{\veps(\cdot,t)},
\]
which is bounded according to \eqref{5.2} combined with Poincaré's inequality and \eqref{5.1}.
\qed
With these preparations at hand, by means of a straightforward subsequence extraction, followed
by applications of the Vitali convergence theorem, Fatou's lemma, Lemma \ref{lem19} and especially Lemma \ref{lem1},
we can now proceed to the construction of a limit pair $(u,\Theta)$ which satisfies (\ref{wu}) and (\ref{wt}).
\begin{lem}\label{lem11}
  There exists $(\eps_j)_{j\in\N}\subset (0,1)$ such that $\eps_j\searrow 0$ as $j\to\infty$, and that with some
  \be{11.1}
	\lbal
	v\in L^\infty_{loc}([0,\infty);L^2(\Om;\R^n)) \cap L^2_{loc}([0,\infty);W_0^{1,2}(\Om;\R^n)), \\[1mm]
	u\in L^\infty_{loc}([0,\infty);W_0^{1,2}(\Om;\R^n)) \qquad \mbox{and} \\[1mm]
	\Theta\in L^\infty_{loc}([0,\infty);L^1(\Om;\R))
		\cap \bigcap_{q\in [1,\frac{n+2}{n})} L^q_{loc}(\bom\times [0,\infty);\R)
		\cap \bigcap_{r\in [1,\frac{n+2}{n+1})} L^r_{loc}([0,\infty);W_0^{1,r}(\Om;\R))
	\ear
  \ee
  fulfilling $\Theta\ge 0$ a.e.~in $\Om\times (0,\infty)$, we have
  \begin{eqnarray}
	& & \veps\to v
	\qquad \mbox{in } L^2_{loc}(\bom\times [0,\infty);\R^n) \mbox{and a.e.~in } \Om\times (0,\infty),
	\label{11.22} \\
	& & \veps\wto v
	\qquad \mbox{in } L^2_{loc}([0,\infty);W_0^{1,2}(\Om;\R^n)),
	\label{11.2} \\
	& & \ueps \to u
	\qquad \mbox{in } L^2_{loc}(\bom\times [0,\infty);\R^n),
	\label{11.3b} \\
	& & \ueps \wto u
	\qquad \mbox{in } L^2_{loc}([0,\infty);W_0^{1,2}(\Om;\R^n))
	\qquad \mbox{and}
	\label{11.3} \\
	& & \Teps\to \Theta
	\qquad \mbox{in } L^1_{loc}(\bom\times [0,\infty);\R) \mbox{and a.e.~in } \Om\times (0,\infty),
	\label{11.4}
  \end{eqnarray}
  as $\eps=\eps_j\searrow 0$.
  These limit functions have the properties that
  \be{11.5}
	u_t=v
	\qquad \mbox{a.e.~in } \Om\times (0,\infty),
  \ee
  that (\ref{wu}) holds for each $\vp\in C_0^\infty(\Om\times [0,\infty);\R^n)$, and that (\ref{wt}) is satisfied for any nonnegative
  $\wh{\vp}\in C_0^\infty(\Om\times [0,\infty);\R)$.
\end{lem}
\proof
  From Lemma \ref{lem5} and Lemma \ref{lem110}, we know that for each $T>0$,
  \bas
	(\veps)_{\eps\in (0,1)}
	\mbox{ is bounded in }
	L^\infty((0,T);L^2(\Om;\R^n))
	\mbox{ and in }
	L^2((0,T);W_0^{1,2}(\Om;\R^n)),
  \eas
  that
  \bas
	(v_{\eps t})_{\eps\in (0,1)}
	\mbox{ is bounded in }
	L^2\big((0,T);(W_0^{2,2}(\Om;\R^n))^\star\big),
  \eas
  that
  \bas
	(\ueps)_{\eps\in (0,1)}
	\mbox{ is bounded in }
	L^\infty((0,T);W_0^{1,2}(\Om;\R^n)),
  \eas
  and that
\bas
	(u_{\eps t})_{\eps\in (0,1)}
	\mbox{ is bounded in }
	L^\infty\big((0,T);(W_0^{2,2}(\Om;\R^n))^\star\big),
  \eas
  while Lemma \ref{lem6} in conjunction with Lemma \ref{lem8}, Lemma \ref{lem9} and Lemma \ref{lem10} says that whenever $T>0$,
  $q\in [1,\frac{n+2}{n})$ and $r\in [1,\frac{n+2}{n+1})$,
  \bas
	(\Teps)_{\eps\in (0,1)}
	\mbox{ is bounded in }
	L^\infty((0,T);L^1(\Om;\R)),
	\mbox{ in }
	L^q(\Om\times (0,T);\R)
	\mbox{ and in }
	L^r((0,T);W_0^{1,r}(\Om;\R)),
  \eas
  and that for each $s>\frac{n+2}{2}$,
  \bas
	(\Theta_{\eps t})_{\eps\in (0,1)}
	\mbox{ is bounded in }
	L^1\big((0,T);(W_0^{2,s}(\Om;\R))^\star\big).
  \eas
  A straightforward extraction procedure involving three applications of an Aubin-Lions lemma (\cite{simon}) combined with
  Vitali's convergence theorem and Fatou's lemma therefore yields $(\eps_j)_{j\in\N}\subset (0,1)$ such that $\eps_j\searrow 0$
  as $j\to\infty$, and that with some functions $v,u$ and $\Theta$ fulfilling (\ref{11.1}) as well as $\Theta\ge 0$ a.e.~in
  $\Om\times (0,\infty)$, the convergence properties in (\ref{11.22})-(\ref{11.4}) hold.\abs
  Now
  according to the second equation in (\ref{0eps}), for each $\phi\in C_0^\infty(\Om\times [0,\infty);\R^n)$ and any $\eps\in (0,1)$
  we have
  \be{11.6}
	- \int_0^\infty\io \ueps\cdot\phi_t
	- \io \ueps(\cdot,0)\cdot\phi(\cdot,0)
	= \eps\int_0^\infty \io \Del\ueps\cdot\phi
	+ \int_0^\infty \io \veps\cdot\phi,
  \ee
  so that since $\eps\Del\ueps\to 0$ in $L^2_{loc}(\bom\times [0,\infty);\R^n)$ as $\eps\searrow 0$ by (\ref{5.6}),
  from (\ref{11.22}), (\ref{ie}) and (\ref{11.3}) we obtain on letting $\eps=\eps_j\searrow 0$ in (\ref{11.6}) that
  \be{11.66}
	- \int_0^\infty \io u\cdot\phi_t
	- \io u_0\cdot\phi(\cdot,0)
	= \int_0^\infty \io v\cdot\phi,
  \ee
  and that thus, in particular, (\ref{11.5}) holds.\abs
  Next, for a verification of (\ref{wu}) we fix $\vp\in C_0^\infty(\Om\times [0,\infty);\R^n)$ and integrate by parts in the first
  equation in (\ref{0eps}) to see that
  \bea{11.7}
	& & \hs{-16mm}
	- \int_0^\infty \io \veps\cdot\vp_t
	- \io v_{0\eps} \cdot \vp(\cdot,0)
	= - \eps \int_0^\infty \io \Del\veps\cdot\Del\vp \nn\\
	& & - \int_0^\infty \io \lan\gamma(\Teps):\nas\veps,\na\vp\ran
	- a \int_0^\infty \io \lan\gamma(\Teps):\nas\ueps,\na\vp\ran,
	\qquad \mbox{for all } \eps\in (0,1),
  \eea
  where by (\ref{11.22}) and (\ref{ie}),
  \be{11.8}
	- \int_0^\infty \io \veps\cdot\vp_t
	- \io v_{0\eps} \cdot \vp(\cdot,0)
	\to
	- \int_0^\infty \io v\cdot\vp_t
	- \io u_{0t} \cdot \vp(\cdot,0)
	\qquad \mbox{as } \eps=\eps_j\searrow 0,
  \ee
  and where
  \be{11.9}
	- \eps \int_0^\infty \io \Del\veps\cdot\Del\vp \to 0
	\qquad \mbox{as } \eps\searrow 0,
  \ee
  because (\ref{5.5}) ensures that also $\eps\Del\veps\to 0$ in $L^2_{loc}(\bom\times [0,\infty);\R^n)$ as $\eps\searrow 0$.
  Moreover, using that $(\gamma(\Teps))_{\eps\in (0,1)}$ is bounded in
  $L^\infty(\Om\times (0,\infty);\R^{n\times n\times n \times n})$ with $\gamma(\Teps)\to\gamma(\Theta)$ a.e.~in
  $\Om\times (0,\infty)$ as $\eps=\eps_j\searrow 0$ by (\ref{11.4}) and the continuity of $\gamma$, from the $L^2$ convergence
  properties of $(\nas v_{\eps j})_{j\in\N}$ and $(\nas u_{\eps_j})_{j\in\N}$ entailed by (\ref{11.2}) and (\ref{11.3}) we infer
  through Lemma \ref{lem19} that
  $\gamma(\Teps):\nas\veps \wto \gamma(\Theta):\nas v$ and
  $\gamma(\Teps):\nas\ueps \wto \gamma(\Theta):\nas u$
  in $L^2_{loc}(\bom\times [0,\infty);\R^{n\times n})$
  as $\eps=\eps_j\searrow 0$.
  Therefore, (\ref{11.7})-(\ref{11.9}) imply that
  \bas
	- \int_0^\infty \io v\cdot\vp_t - \io u_{0t}\cdot\vp(\cdot,0)
	= - \int_0^\infty \io \lan\gamma(\Theta):\nas v,\na\vp\ran
	- a \int_0^\infty \io \lan\gamma(\Theta):\nas u,\na\vp\ran,
  \eas
  whence (\ref{wu}) follows upon observing that $\nas v=\nas u_t$ by (\ref{11.5}), and that according to (\ref{11.66}) when
  applied to $\phi:=\vp_t$,
  \bas
	- \int_0^\infty \io v\cdot\vp_t
	= \int_0^\infty \io u\cdot\vp_{tt}
	+ \io u_0 \cdot \vp_t(\cdot,0).
  \eas
  Finally, for arbitrary nonnegative $\wh{\vp}\in C_0^\infty(\Om\times [0,\infty);\R)$ and $\eps\in (0,1)$ we may use the third
  equation in (\ref{0eps}) to see that
  \be{11.10}
	\int_0^\infty \io \lan\Gamma(\Teps):\nas\veps,\nas\veps\ran \wh{\vp}
	= - \int_0^\infty \io \Teps \wh{\vp}_t
	- \io \Theta_{0\eps} \wh{\vp}(\cdot,0)
	- D \int_0^\infty \io \Teps \Del\wh{\vp},
  \ee
  and note that here, due to (\ref{11.4}) and (\ref{ie}),
  \bea{11.11}
	& & \hs{-36mm}
	- \int_0^\infty \io \Teps \wh{\vp}_t
	- \io \Theta_{0\eps} \wh{\vp}(\cdot,0)
	- D \int_0^\infty \io \Teps \Del\wh{\vp} \nn\\
	&\to& - \int_0^\infty \io \Theta \wh{\vp}_t
	- \io \Theta_0 \wh{\vp}(\cdot,0)
	- D \int_0^\infty \io \Theta \Del\wh{\vp}
  \eea
  as $\eps=\eps_j\searrow 0$.
  Since in view of (\ref{Gamma_symm}), (\ref{Gamma_pos}), (\ref{11.4}) and (\ref{11.2}) we may employ Lemma \ref{lem1} to infer
  thanks to the boundedness of $\supp \wh{\vp}$ that
  \bas
	\int_0^\infty \io \lan\Gamma(\Theta):\nas v,\nas v\ran \wh{\vp}
	\le \liminf_{\eps=\eps_j\searrow 0}
	\int_0^\infty \io \lan\Gamma(\Teps):\nas\veps,\nas\veps\ran \wh{\vp},
  \eas
  from (\ref{11.10}) and (\ref{11.11}) we obtain that indeed the inequality in (\ref{wt}) holds,
  whereby the proof is completed.
\qed
The most crucial step in our analysis, however, can be found in the following argument	
which now makes full use of the results from Section \ref{sect_lsc} to derive (\ref{wF}) after an appropriate
rearrangement of the key contributions to the corresponding dissipation rate (see (\ref{D}))
as essentially quadratic expressions  (cf.~(\ref{12.14})).
\begin{lem}\label{lem12}
  Let $u$ and $\Theta$ be as in Lemma \ref{lem11}. Then there exist $\kappa>0,\lam>0$ and $\mu>0$ such that (\ref{wF})
  is valid for each nonnegative $\psi\in C_0^\infty(\Om;\R)$ and any nonincreasing
  $\zeta\in C_0^\infty([0,\infty);\R)$.
\end{lem}
\proof
  We let $K_\gamma>0$ be as in (\ref{gamma_pos}), and taking $\eta_1(\cdot)$ and $\eta_2(\cdot)$ from Lemma \ref{lem102}
  and Lemma \ref{lem103}, we fix $\kappa>0$ in such a way that
  \be{12.1}
	\frac{\kappa}{a}<\eta_2(\gamma).
  \ee
  We thereupon choose $\lam>0$ small enough fulfilling
  \be{12.2}
	\frac{\kappa}{a} \eta_1(\Gamma) >\lam
  \ee
  and pick some $\mu>0$ suitably large such that
  \be{12.3}
	\frac{\kappa(a+2\mu)}{4} \eta_1(\gamma)> \frac{a^2}{4}.
  \ee
  To verify that then (\ref{wF}) holds for any nonnegative $\psi\in C_0^\infty(\Om;\R)$ and any
  $\zeta\in C_0^\infty([0,\infty);\R)$ with $\zeta_t\le 0$, letting
  \be{12.33}
	\Feps\equiv \Feps^{(\kappa,\lam,\psi)}(x,t)
	:= \Big(\frac{1}{2} |\veps|^2 + \frac{\kappa}{2} |\na\ueps|^2 + \lam\Teps \Big) \psi
  \ee
  we go back to Lemma \ref{lem2}, Lemma \ref{lem3} and Lemma \ref{lem4} and thereby see that since
  $\supp \zeta$ is bounded and $\psi|_{\pO}=0$,
  \bea{12.4}
	& & \hs{-24mm}
	\ - \ \zeta(0) \io \Feps(\cdot,0)
	= \int_0^\infty \io \pa_t \big\{ \Feps \zeta(t) e^{-\mu t} \big\} \nn\\
	&=& \int_0^\infty \io \Big\{ \frac{1}{2} (|\veps|^2)_t \psi
	+ \frac{\kappa}{2} (|\na\ueps|^2)_t \psi
	+ \lam \Theta_{\eps t} ψ \Big\} \zeta(t) e^{-\mu t} \nn\\
	& & + \int_0^\infty \io \Big\{ \frac{1}{2} |\veps|^2 \psi + \frac{\kappa}{2} |\na\ueps|^2 \psi + \lam\Teps\psi\Big\}
		\big\{ \zeta_t(t)-\mu\zeta(t) \big\} e^{-\mu t} \nn\\
	&=& \int_0^\infty \io \Big\{ - \lan\gamma(\Teps):\nas\veps,\nas\veps\ran\psi
	-\eps |\Del\veps|^2 \psi \nn\\
	& & \hs{17mm}
	- a\lan\gamma(\Teps):\nas\ueps,\nas\veps\ran\psi
	+ \lan\gamma(\Teps):\nas\veps,\veps\mult \na\psi\ran \nn\\
	& & \hs{17mm}
	+ a\lan\gamma(\Teps):\nas\ueps,\veps\mult\na\psi\ran \nn\\
	& & \hs{17mm}
	-2\eps (\na\veps\cdot\na\psi)\Del\veps
	- \eps(\veps\cdot\Del\veps)\Del\psi \Big\} \zeta(t) e^{-\mu t} \nn\\
	& & + \kappa \int_0^\infty \io \Big\{
	-\eps |\Del\ueps|^2 \psi
	+ \lan\na\ueps,\na\veps\ran\psi
	- \eps (\na\ueps\cdot\na\psi)\cdot\Del\ueps
	\Big\} \zeta(t) e^{-\mu t} \nn\\
	& & + \lam \int_0^\infty \io \Big\{
	D\Teps \Del\psi
	+ \lan\Gamma(\Teps):\nas\veps,\nas\veps\ran \psi
	\Big\} \zeta(t) e^{-\mu t} \nn\\
	& & + \int_0^\infty \io \Big\{
	\frac{1}{2} |\veps|^2 \psi + \frac{\kappa}{2} |\na\ueps|^2 \psi + \lam\Teps \psi
	\Big\} \big\{ \zeta_t(t) - \mu\zeta(t)\big\} e^{-\mu t} \nn\\
	&=& - \int_0^\infty \io \lan\gamma(\Teps):\nas\veps,\nas\veps\ran \zeta(t) e^{-\mu t} \psi
	- a \int_0^\infty \io \lan\gamma(\Teps):\nas\ueps,\nas\veps\ran \zeta(t) e^{-\mu t} \psi \nn\\
	& & - \frac{\kappa\mu}{2} \int_0^\infty\io |\na\ueps|^2 \zeta(t) e^{-\mu t} \psi
	+ \kappa \int_0^\infty \io \lan\na\ueps,\na\veps\ran \zeta(t) e^{-\mu t} \psi \nn\\
	& & + \lam \int_0^\infty \io \lan\Gamma(\Teps):\nas\veps,\nas\veps\ran \zeta(t) e^{-\mu t} \psi \nn\\
	& & - J^{(1)}(\eps)
	+ J^{(2)}(\eps)
	\qquad \mbox{for all } \eps\in (0,1),
  \eea
  where
  \bea{12.5}
	J^{(1)}(\eps)
	&:=& \int_0^\infty \io \Big( \frac{1}{2} |\veps|^2 + \lam\Teps\Big) (\mu\zeta(t)-\zeta_t(t)) e^{-\mu t} \psi
	- \frac{\kappa}{2} \int_0^\infty \io |\na\ueps|^2 \zeta_t(t) e^{-\mu t} \psi \nn\\
	& & + \eps \int_0^\infty \io |Δ\veps|^2 \zeta(t) e^{-\mu t} \psi
	+ \kappa\eps \int_0^\infty \io |\Del\ueps|^2 \zeta(t) e^{-\mu t} \psi,
	\qquad \eps\in (0,1),
  \eea
  and
  \bea{12.6}
	J^{(2)}(\eps)
	&:=& \int_0^\infty \io \lan\gamma(\Teps):\nas\veps,\veps\mult\na\psi\ran \zeta(t) e^{-\mu t}
	+ a \int_0^\infty \io \lan \gamma(\Teps):\nas\ueps,\veps\mult\na\psi\ran \zeta(t) e^{-\mu t} \nn\\
	& & + \lam D \int_0^\infty \io \Teps \zeta(t) e^{-\mu t} \Del\psi \nn\\
	& & - 2\eps \int_0^\infty \io \big\{ (\na\veps\cdot\na\psi)\cdot\Del\veps\big\} \zeta(t) e^{-\mu t}
	- \eps \int_0^\infty \io (\veps\cdot\Del\veps) \zeta(t) e^{-\mu t} \Del\psi \nn\\
	& & - \kappa\eps \int_0^\infty \io \big\{ (\na\ueps\cdot\na\psi)\cdot\Del\ueps\big\} \zeta(t) e^{-\mu t},
	\qquad \eps\in (0,1).
  \eea
  Here since Lemma \ref{lem19} along with (\ref{11.2}), (\ref{11.4})
  and the boundedness of $\gamma$ again ensures that with $(\eps_j)_{j\in\N}$ taken from
  Lemma \ref{lem11} we have
  $\gamma(\Teps):\nas\veps \wto \gamma(\Theta):\nas v$
  in $L^2_{loc}(\bom\times [0,\infty);\R^{n\times n})$ as $\eps=\eps_j\searrow 0$, the strong convergence property in
  (\ref{11.22}) implies that
  \be{12.7}
	\int_0^\infty \io \lan\gamma(\Teps):\nas\veps,\veps\mult\na\psi\ran \zeta(t) e^{-\mu t}
	\to \int_0^\infty \io \lan\gamma(\Theta):\nas v,v\mult\na\psi\ran \zeta(t) e^{-\mu t}
	\qquad \mbox{as } \eps=\eps_j\searrow 0,
  \ee
  because $T:=\sup\{t>0 \ | \ \zeta(t)\ne 0\}$ is finite.
  Similarly, (\ref{11.3}), (\ref{11.4}) and (\ref{11.22}) guarantee that
  \be{12.8}
	a \int_0^\infty \io \lan \gamma(\Teps):\nas\ueps,\veps\mult\na\psi\ran \zeta(t) e^{-\mu t}
	\to a \int_0^\infty \io \lan\gamma(\Theta):\nas u,v\mult\na\psi\ran \zeta(t) e^{-\mu t}
	\qquad \mbox{as } \eps=\eps_j\searrow 0,
  \ee
  while (\ref{11.4}) and the compactness of $\supp\Del\psi\subset\Om$ entail that
  \be{12.9}
	\lam D \int_0^\infty \io \Teps\zeta(t) e^{-\mu t} \Del\psi
	\to \lam D \int_0^\infty \io \Theta \zeta(t) e^{-\mu t} \Del\psi
	\qquad \mbox{as } \eps=\eps_j\searrow 0.
  \ee
  To see that the last three summands vanish in the limit $\eps\searrow 0$, we abbreviate $c_1:=\|\zeta\|_{L^\infty([0,\infty))}$,
  $c_2:=\|\na\psi\|_{L^\infty(\Om)}$ and $c_3:=\|\Del\psi\|_{L^\infty(\Om)}$, and use the Cauchy-Schwarz inequality in estimating
  \bas
	\bigg| -2\eps \int_0^\infty \io \big\{ (\na\veps\cdot\na\psi)\cdot\Del\veps\big\} \zeta(t) e^{-\mu t} \bigg|
	&\le& 2c_1 c_2 \eps \int_0^T \io |\na\veps| \, |\Del\veps| \\
	&\le& 2c_1 c_2 \eps^\frac{1}{2} \bigg\{ \eps \int_0^T \io |\Del\veps|^2 \bigg\}^\frac{1}{2}
		\bigg\{ \int_0^T \io |\na\veps|^2 \bigg\}^\frac{1}{2}
  \eas
  and
  \bas
	\bigg| -\eps \int_0^\infty \io (\veps\cdot\Del\veps) \zeta(t) e^{-\mu t} \Del\psi \bigg|
	&\le& c_1 c_3 \eps \int_0^T \io |\veps| \, |\Del\veps| \\
	&\le& c_1 c_3 \eps^\frac{1}{} \bigg\{ \eps\int_0^T \io |\Del\veps|^2 \bigg\}^\frac{1}{2}
		\bigg\{ \int_0^T \io |\veps|^2 \bigg\}^\frac{1}{2}
  \eas
  as well as
  \bas
	\bigg| -\kappa\eps \int_0^\infty \io \big\{ (\na\ueps\cdot\na\psi)\cdot\Del\ueps\big\} \zeta(t) e^{-\mu t} \bigg|
	&\le& c_1 c_2\kappa \eps \int_0^T \io |\na\ueps| \, |\Del\ueps| \\
	&\le& c_1 c_2 \kappa \eps^\frac{1}{2} \bigg\{ \eps\int_0^T \io |\Del\ueps|^2 \bigg\}^\frac{1}{2}
		\bigg\{ \int_0^T \io |\na\ueps|^2 \bigg\}^\frac{1}{2}
  \eas
  for all $\eps\in (0,1)$.
  In consequence, from (\ref{5.5}), (\ref{5.4}), (\ref{5.1}), (\ref{5.6}) and (\ref{5.2}) we thus infer that, indeed,
  \bas
	- 2\eps \int_0^\infty \io \big\{ (\na\veps\cdot\na\psi)\cdot\Del\veps\big\} \zeta(t) e^{-\mu t}
	- \eps \int_0^\infty \io (\veps\cdot\Del\veps) \zeta(t) e^{-\mu t} \Del\psi \nn\\
	- \kappa\eps \int_0^\infty \io \big\{ (\na\ueps\cdot\na\psi)\cdot\Del\ueps\big\} \zeta(t) e^{-\mu t}
	\to 0
	\qquad \mbox{as } \eps\searrow 0,
  \eas
  whence (\ref{12.6})-(\ref{12.9}) imply that
  \bea{12.10}
	J^{(2)}(\eps)
	&\to& \int_0^\infty \io \lan\gamma(\Theta):\nas v,v\mult\na\psi\ran \zeta(t) e^{-\mu t}
	+ a \int_0^\infty \io \lan\gamma(\Theta):\nas u,v\mult\na\psi\ran \zeta(t) e^{-\mu t} \nn\\
	& & + \lam D \int_0^\infty \io \Theta\zeta(t) e^{-\mu t} \Del\psi
	\qquad \mbox{as } \eps=\eps_j\searrow 0.
  \eea
  Next addressing (\ref{12.5}), we may rely on the fact that both $\zeta$ and $\psi$ are nonnegative to simply estimate
  \be{12.11}
	\eps \int_0^\infty \io |\Del\veps|^2 \zeta(t) e^{-\mu t} \psi
	+ \kappa\eps \int_0^\infty \io |\Del\ueps|^2 \zeta(t) e^{-\mu t} \psi
	\ge 0
	\qquad \mbox{for all } \eps\in (0,1),
  \ee
  and use the pointwise approximation features contained in (\ref{11.22}) and (\ref{11.4}) to see that since also $-\zeta_t$ is
  nonnegative, Fatou's lemma implies that
  \be{12.12}
	\int_0^\infty \io \Big(\frac{1}{2} |v|^2 + \lam\Theta\Big) \big(\mu\zeta(t)-\zeta_t(t)\big) e^{-\mu t} \psi
	\le \liminf_{\eps=\eps_j\searrow 0} \int_0^\infty \io
		\Big(\frac{1}{2} |\veps|^2 + \lam\Teps\Big)\big(\mu\zeta(t)-\zeta_t(t)\big) e^{-\mu t} \psi.
  \ee
  Since in view of (\ref{11.3}) the lower semicontinuity of the norm in $L^2(\Om\times (0,T);\R^{n\times n})$ ensures that
  \bas
	-\frac{\kappa}{2} \int_0^\infty \io |\na u|^2 \zeta_t(t) e^{-\mu t} \psi
	&=& \frac{\kappa}{2} \int_0^T \io \Big| \sqrt{|\zeta_t(t)| e^{-\mu t} \psi} \na u \Big|^2 \\
	&\le& \frac{\kappa}{2} \liminf_{\eps=\eps_j\searrow 0}
		\int_0^T \io \Big| \sqrt{|\zeta_t(t)| e^{-\mu t} \psi} \na \ueps \Big|^2 \\
	&=& \liminf_{\eps=\eps_j\searrow 0} \bigg\{
		-\frac{\kappa}{2} \int_0^\infty \io |\na\ueps|^2 \zeta_t(t) e^{-\mu t} \psi \bigg\},
  \eas
  a combination of (\ref{12.11}) with (\ref{12.12}) shows that
  \bea{12.13}
	& & \hs{-20mm}
	\int_0^\infty \io \Big(\frac{1}{2}|v|^2 + \lam\Theta\Big)\big(\mu\zeta(t)-\zeta_t(t)\big) e^{-\mu t} \psi
	- \frac{\kappa}{2} \int_0^\infty \io |\na u|^2 \zeta_t(t) e^{-\mu t} \psi \nn\\
	&\le& \liminf_{\eps=\eps_j\searrow 0} J^{(1)}(\eps).
  \eea

  Henceforth directing our attention toward the crucial first five summands on the right of (\ref{12.4}), we create quadratic
  expression therein by observing that thanks to the symmetriy property (\ref{gamma_symm}),
  \bas
	\lan\gamma(\Teps):\nas (\veps+\mbox{$\frac{a}{2}$}\ueps),\nas (\veps+\mbox{$\frac{a}{2}$}\ueps) \ran
	&=& \lan\gamma(\Teps):\nas\veps,\nas\veps\ran
	+ \frac{a}{2} \lan\gamma(\Teps):\nas\ueps,\nas\veps\ran \\
	& & + \frac{a}{2} \lan\gamma(\Teps):\nas\veps,\nas\ueps\ran
	+ \frac{a^2}{4} \lan\gamma(\Teps):\nas\ueps,\nas\ueps\ran \\
	&=& \lan\gamma(\Teps):\nas\veps,\nas\veps\ran
	+ a \lan\gamma(\Teps):\nas\ueps,\nas\veps\ran \\
	& & + \frac{a^2}{4} \lan\gamma(\Teps):\nas\ueps,\nas\ueps\ran
  \eas
  and hence
  \bas
	& & \hs{-20mm}
	\lan\gamma(\Teps):\nas\veps,\nas\veps\ran
	+ a \lan\gamma(\Teps):\nas\ueps,\nas\veps\ran \nn\\
	&=&
	\lan\gamma(\Teps):\nas (\veps+\mbox{$\frac{a}{2}$}\ueps),\nas (\veps+\mbox{$\frac{a}{2}$}\ueps) \ran
	- \frac{a^2}{4} \lan\gamma(\Teps):\nas\ueps,\nas\ueps\ran
  \eas
  in $\Om\times (0,\infty)$ for all $\eps\in (0,1)$, and that similarly,
  \bas
	- \kappa\lan\na\ueps,\na\veps\ran
	=
	- \frac{\kappa}{a} \big|\na\big(\veps+\mbox{$\frac{a}{2}$} \ueps\big)\big|^2
	+ \frac{\kappa}{a} |\na\veps|^2
	+ \frac{\kappa a}{4} |\na\ueps|^2
  \eas
  in $\Om\times (0,\infty)$ for all $\eps\in (0,1)$.
  Therefore, indeed,
  \bea{12.14}
	& & \hs{-12mm}
	\int_0^\infty \io \lan\gamma(\Teps):\nas\veps,\nas\veps\ran \zeta(t) e^{-\mu t} \psi
	+ a \int_0^\infty \io \lan\gamma(\Teps):\nas\ueps,\nas\veps\ran \zeta(t) e^{-\mu t} \psi \nn\\
	& & \hs{-12mm}
	+ \frac{\kappa\mu}{2} \int_0^\infty\io |\na\ueps|^2 \zeta(t) e^{-\mu t} \psi
	- \kappa \int_0^\infty \io \lan\na\ueps,\na\veps\ran \zeta(t) e^{-\mu t} \psi \nn\\
	& & \hs{-12mm}
	- \lam \int_0^\infty \io \lan\Gamma(\Teps):\nas\veps,\nas\veps\ran \zeta(t) e^{-\mu t} \psi \nn\\
	&=& \int_0^\infty \io \lan\gamma(\Teps):\nas (\veps+\mbox{$\frac{a}{2}$}\ueps),\nas (\veps+\mbox{$\frac{a}{2}$}\ueps) \ran
		\zeta(t) e^{-\mu t} \psi
	- \frac{\kappa}{a} \int_0^\infty \io \big|\na\big(\veps+\mbox{$\frac{a}{2}$} \ueps\big)\big|^2 \zeta(t) e^{-\mu t} \psi \nn\\
	& & + \frac{\kappa}{a} \int_0^\infty \io |\na\veps|^2 \zeta(t) e^{-\mu t} \psi
	- \lam \int_0^\infty \io \lan\Gamma(\Teps):\nas\veps,\nas\veps\ran \zeta(t) e^{-\mu t} \psi \nn\\
	& & + \frac{\kappa(a+2\mu)}{4} \int_0^\infty \io |\na\ueps|^2 \zeta(t) e^{-\mu t} \psi
	- \frac{a^2}{4} \int_0^\infty \io \lan\gamma(\Teps):\nas\ueps,\nas\ueps\ran \zeta(t) e^{-\mu t} \psi
  \eea
  for all $\eps\in (0,1)$.
  Now (\ref{12.1}) together with (\ref{gamma_pos}) and (\ref{11.22})-(\ref{11.4}) enables us to see upon an application of
  Lemma \ref{lem103} to $\beta:=\gamma$, to $w_j:=\veps+\frac{a}{2}\ueps$
  and $z_j:=\Teps$ for $\eps=\eps_j$, and to
  $\vp(x,t):=\zeta(t) e^{-\mu t} \psi(x)$ for $(x,t)\in \bom\times [0,T]$ that
  \bas
	& & \hs{-20mm}
	\int_0^\infty \io \lan\gamma(\Theta):\nas (v+\mbox{$\frac{a}{2}$}u)\nas (v+\mbox{$\frac{a}{2}$}u)\ran
		\zeta(t) e^{-\mu t} \psi
	- \frac{\kappa}{a} \int_0^\infty \io \big|\na\big(v+\mbox{$\frac{a}{2}$}u\big)\big|^2 \zeta(t) e^{-\mu t} \psi \\
	& & \le \liminf_{\eps=\eps_j\searrow 0} \bigg\{
	\int_0^\infty \io \lan\gamma(\Teps):\nas (\veps+\mbox{$\frac{a}{2}$}\ueps),\nas(\veps+\mbox{$\frac{a}{2}$}\ueps)\ran
		\zeta(t) e^{-\mu t} \psi \\
	& & \hs{20mm}
	- \frac{\kappa}{a} \int_0^\infty \io \big|\na\big(\veps+\mbox{$\frac{a}{2}$}\ueps\big)\big|^2 \zeta(t) e^{-\mu t} \psi
	\bigg\},
  \eas
  while similarly from (\ref{12.2}) and Lemma \ref{lem102} we obtain that
  \bas
	& & \hs{-20mm}
	\frac{\kappa}{a} \int_0^\infty \io |\na v|^2 \zeta(t) e^{-\mu t} \psi
	- \lam \int_0^\infty \io \lan\Gamma(\Theta):\nas v,\nas v\ran \zeta(t) e^{-\mu t} \psi \\
	& & \le \liminf_{\eps=\eps_j\searrow 0} \bigg\{
	\frac{\kappa}{a} \int_0^\infty \io |\na\veps|^2 \zeta(t) e^{-\mu t} \psi
	- \lam \int_0^\infty \io \lan\Gamma(\Teps):\nas\veps,\nas\veps\ran \zeta(t) e^{-\mu t} \psi
	\bigg\},
  \eas
  and while (\ref{12.3}) in view of Lemma \ref{lem102} ensures that
  \bas
	& & \hs{-12mm}
	\frac{\kappa(a+2\mu)}{4} \int_0^\infty \io |\na u|^2 \zeta(t) e^{-\mu t} \psi
	- \frac{a^2}{4} \int_0^\infty \io \lan\gamma(\Theta):\nas u,\nas u\ran \zeta(t) e^{-\mu t} \psi \\
	& & \le \liminf_{\eps=\eps_j\searrow 0} \bigg\{
	\frac{\kappa(a+2\mu)}{4} \int_0^\infty \io |\na\ueps|^2 \zeta(t) e^{-\mu t} \psi
	- \frac{a^2}{4} \int_0^\infty \io \lan\gamma(\Teps):\nas\ueps,\nas\ueps\ran \zeta(t) e^{-\mu t} \psi
	\bigg\}.
  \eas
  Since clearly
  \bas
	\zeta(0) \io \Feps(\cdot,0)
	&=& \zeta(0) \io \Big( \frac{1}{2} |v_{0\eps}|^2 + \frac{\kappa}{2} |\na u_{0\eps}|^2 + \lam \Theta_{0\eps} \Big) \psi \\
	&\to& \zeta(0) \io \Big( \frac{1}{2} |u_{0t}|^2 + \frac{\kappa}{2} |\na u_0|^2 + \lam\Theta_0\Big) \psi
	\qquad \mbox{as } \eps\searrow 0
  \eas
  by (\ref{ie}), in conjunction with (\ref{12.10}) and (\ref{12.13}) this shows that (\ref{12.4}) and (\ref{12.14}) imply
  the inequality
  \bea{12.99}
	& & \hs{-16mm}
	\int_0^\infty \io \lan\gamma(\Theta):\nas(v+\mbox{$\frac{a}{2}$}u)\nas (v+\mbox{$\frac{a}{2}$}u)\ran \zeta(t) e^{-\mu t} \psi
	- \frac{\kappa}{a} \int_0^\infty \io \big|\na\big(v+\mbox{$\frac{a}{2}$}u\big)\big|^2 \zeta(t) e^{-\mu t} \psi \nn\\
	& & \hs{-16mm}
	+ \frac{\kappa}{a} \int_0^\infty \io |\na v|^2 \zeta(t) e^{-\mu t} \psi
	- \lam \int_0^\infty \io \lan\Gamma(\Theta):\nas v,\nas v\ran \zeta(t) e^{-\mu t} \psi \nn\\
	& & \hs{-16mm}
	\frac{\kappa(a+2\mu)}{4} \int_0^\infty \io |\na u|^2 \zeta(t) e^{-\mu t} \psi
	- \frac{a^2}{4} \int_0^\infty \io \lan\gamma(\Theta):\nas u,\nas u\ran \zeta(t) e^{-\mu t} \psi \nn\\
	& & \hs{-16mm}
	+ \int_0^\infty \io \Big(\frac{1}{2} |v|^2 + \lam\Theta\Big)\big(\mu\zeta(t)-\zeta_t(t)\big) e^{-\mu t} \psi
	- \frac{\kappa}{2} \int_0^\infty \io |\na u|^2 \zeta_t(t) e^{-\mu t} \psi \nn\\
	&\le& \liminf_{\eps=\eps_j\searrow 0} \bigg\{
	\int_0^\infty \io \lan\gamma(\Teps):\nas\veps,\nas\veps\ran \zeta(t) e^{-\mu t} \psi
	+ a \int_0^\infty \io \lan\gamma(\Teps):\nas\ueps,\nas\veps\ran \zeta(t) e^{-\mu t} \psi \nn\\
	& & \hs{15mm}
	+ \frac{\kappa\mu}{2} |\na\ueps|^2 \zeta(t) e^{-\mu t} \psi
	- \kappa \int_0^\infty \io \lan\na\ueps,\na\veps\ran \zeta(t) e^{-\mu t} \psi \nn\\
	& & \hs{15mm}
	- \lam \int_0^\infty \io \lan\Gamma(\Teps):\nas\veps,\nas\veps\ran \zeta(t) e^{-\mu t} \psi \nn\\
	& & \hs{15mm}
	+ J^{(1)}(\eps) \bigg\} \nn\\
	&=& \liminf_{\eps=\eps_j\searrow 0} \bigg\{
	\zeta(0) \io \Feps(\cdot,0) + J^{(2)}(\eps)\bigg\} \nn\\
	&=& \zeta(0) \io \Big(\frac{1}{2} |u_{0t}|^2 + \frac{\kappa}{2} |\na u_0|^2 + \lam\Theta_0\Big) \psi \nn\\
	& & + \int_0^\infty \io \lan\gamma(\Theta):\nas v,v\mult\na\psi\ran \zeta(t) e^{-\mu t}
	+ a \int_0^\infty \io \lan\gamma(\Theta):\nas u,v\mult\na\psi\ran \zeta(t) e^{-\mu t}  \nn\\
	& & + \lam D \int_0^\infty \io \Theta\zeta(t) e^{-\mu t} \Del\psi.
  \eea
  In line with (\ref{F0}), (\ref{D}), (\ref{F}) and (\ref{R}), rewriting the left-hand side herein in the style of (\ref{12.14})
  leads to (\ref{wF}), because $u_t=v$ a.e.~in $\Om\times (0,\infty)$ by (\ref{11.5}).
\qed
The proof of our main result has thereby actually been completed:\abs
\proofc of Theorem \ref{theo13}. \quad
  Taking $(u,\Theta)$ from Lemma \ref{lem11}, we only need to combine the latter with Lemma \ref{lem12}.
\qed
\section{Numerical experiments}
This section intends to conduct some numerical experiments to illustrate some possible influences
of temperature-dependent material parameters on the results of acoustic processes.
Specifically, the coupled thermal and mechanical behavior of a one-dimensional acoustic resonator is considered using a finite-difference time-domain (FDTD) method~\cite{Yee1966,Claes2024},
and to facilitate a connection to corresponding literature we shall here return to the original variables by
considering the evolution system
\begin{align}
	\rho u_{tt} &= \tau C(\Theta) u_{xxt} + C(\Theta) u_{xx}, \nonumber \\
	c\rho \Theta_t &= \lam \Theta_{xx} + \tau C(\Theta) u_{xt}^2,
    \label{num}
\end{align}
in which all considered quantities are scalar.\abs
The setup roughly approximates an actively cooled (\(\Theta = 0\) at the boundaries), mechanically clamped (\(u = u_t = 0\) at the boundaries) layer that is driven to oscillate at its resonance frequency.
The viscous wave equation and the heat equation in (\ref{num})
are solved with parameter values listed in Table~\ref{tab:parameters}.
The elasticity \(C\) of the material is assumed to be temperature dependent, using two different models.
To explore temperature dependent behavior that is akin to the observations recorded in Figure~\ref{fig:elastic_measurement},
a power-law in temperature is used,
\begin{equation}
    C = C_0 \cdot (1 + k \Theta ^ p),
    \label{equ:temp_dependence_pow}
\end{equation}
neglecting here the boundedness requirements made in Theorem~\ref{theo13}.
For the analysis of convergent behavior, an exponential law is presupposed,
\begin{equation}
    C = C_0 \cdot (\al + (1 - \al) e^{-b \Theta}) ,
    \label{equ:temp_dependence_exp}
\end{equation}
which yields \(C = C_0\) for \(\Theta = 0\) and converges to \(C = \al \cdot C_0\) in the limit of large positive values of \(\Theta \).\abs
The thickness of the resonator is 1\,mm, resulting in a resonance frequency (at \(\Theta = 0\)) of 2\,MHz.
A continuous sinusoidal signal at this frequency, implemented as a time-depended Dirichlet boundary condition for the velocity \(u_t\) at the centre of the resonator, is applied to excite the system.
Due to the high mechanical frequency and comparatively slow thermal processes, the simulation requires large number of steps (\(5 \cdot 10^5\)) in time-domain to show significant effects.
The spatial and temporal resolution are chosen to over-satisfy the conditions of stability for the solution of wave equations~\cite{Remis2000} by a factor of 2.5 (\(\Delta x / \Delta t = 2.5 \cdot c_\mathrm{ph}\)).

\begin{figure}
    \centering
    \input{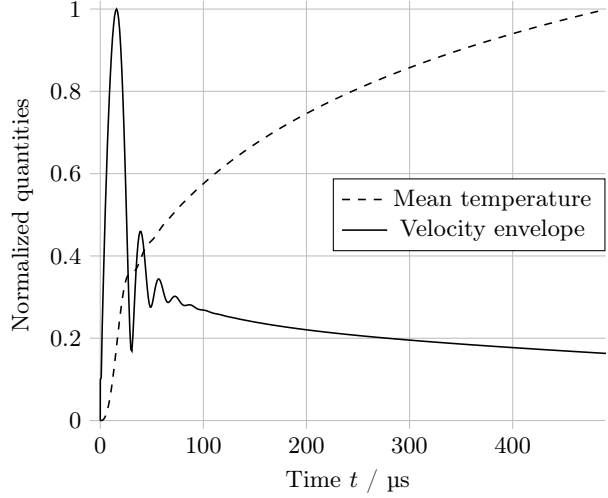}
    \caption{Typical result for the mean temperature and for the envelope of the mechanical oscillation from a coupled thermo-acoustic simulations of a one-dimensional resonator. Parameters for the temperature dependence power-law of the elasticity are \(k = 10^7\) and \(p = 1\).}\label{fig:temp_envelope}
\end{figure}

\begin{table}
    \caption{Material parameters for the numerical study. Values approximate the mechanical and thermal properties of a piezoelectric ceramic~\cite{Feldmann2021,PICeramic2023}.}\label{tab:parameters}
    \centering
    \begin{tabular}{c c c}
        Parameter & Value & Unit \\
        \hline
        \(C_0\) & 124.8 & GPa \\
        \(\rho \) & 7800 & kg\,m\(^{-3}\) \\
        \(\tau \) & 1 & ns \\
        \(c\) & 350 & J\,K\(^{-1}\) \,kg\(^{-1}\) \\
        \(\lambda \) & 1.1 & W\,m\(^{-1}\) \,K\(^{-1}\) \\
        \hline
    \end{tabular}
\end{table}
To analyse the mechanical and thermal behavior of the resonator, the temperature and the velocity \(u_t\) are observed.
For a clearer depiction, the normalized mean temperature and the envelope of the velocity signal are shown.
The results presented in \autoref{fig:temp_envelope} constitute an archetypical behavior observed for many parameter value  tuples \((k, p)\) and \((\al, b)\):
Initially, the velocity of the oscillation increases rapidly because of resonance.
The mean temperature of the material initially shows superlinear increase due to the mechanical losses, causing the elasticity to change and thus a shift in the resonance frequency of the oscillator.
Because the excitation remains at 2\,MHz, the system is no longer excited in resonance, leading to a reduced oscillation velocity after a short period of beating.
The trend of the mean temperature increase thus changes to a sublinear behavior.
The overall behavior of the mean temperature and the velocity envelope is similar for large ranges of \(k\) and \(p\) of the power law as well as \(\al\) and \(b\) for the exponential expression, with only quantitative variations, e.g.\ when the velocity envelope begins to decrease.
However, there are a number of configurations, especially for large values of \(k\), where the simulation destabilises and the temperature field values overflow.
Similar behavior can be observed when using the convergent, exponential expression to model the temperature dependence, however, overflow can only be brought about if the increase in elasticity occurs sudden and steep at the start of the simulation, e.g.\ both parameters \(\al\) and \(b\) need to have high values.
It is still to be shown if this behavior results from numerical problems with the finite-difference scheme or is inherent to the system of equations.

\begin{figure}
    \centering
    \begin{tikzpicture}
    \footnotesize

    \begin{axis}[
        axis line style={draw=none},
        width=\textwidth/3,
        title={a: \(k = 10^2;\ p = 0.5\)},
        name=plot0,
        ylabel={Position \(x\) / mm},
        xmin=0, xmax=499,
        ymin=0, ymax=1,
        xticklabels={,,},
        title style={yshift=-7},
        tick align=outside,
        tick pos=left,
    ]
        \addplot graphics [xmin=0, xmax=499, ymin=0, ymax=1] {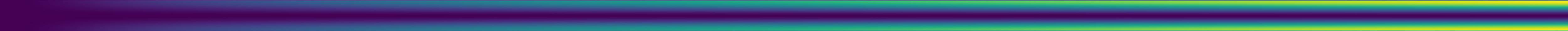};
    \end{axis}

    \begin{axis}[
        axis line style={draw=none},
        width=\textwidth/3,
        title={b: \(k = 10^4;\ p = 0.5\)},
        name=plot1,
        at=(plot0.right of south east), anchor=left of south west,
        xmin=0, xmax=499,
        ymin=0, ymax=1,
        yticklabels={,,},
        xticklabels={,,},
        title style={yshift=-7},
        tick align=outside,
        tick pos=left,
    ]
        \addplot graphics [xmin=0, xmax=499, ymin=0, ymax=1] {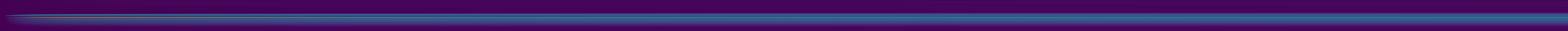};
    \end{axis}

    \begin{axis}[
        axis line style={draw=none},
        width=\textwidth/3,
        title={c: \(k = 10^7;\ p = 1\)},
        name=plot2,
        at=(plot1.right of south east), anchor=left of south west,
        xmin=0, xmax=499,
        ymin=0, ymax=1,
        yticklabels={,,},
        xticklabels={,,},
        title style={yshift=-7},
        tick align=outside,
        tick pos=left,
        point meta min=0,
        point meta max=1,
        colorbar,
        colormap/viridis,
        colorbar style={
            axis line style={draw=none},
            ymin=0, ymax=1,
            ylabel={Norm. Temperature},
            tick align=outside,
            tick pos=right,
            at={(1.05,0.5)},
            anchor=west,
            width=0.04*\pgfkeysvalueof{/pgfplots/parent axis width},
        },
    ]
        \addplot graphics [xmin=0, xmax=499, ymin=0, ymax=1] {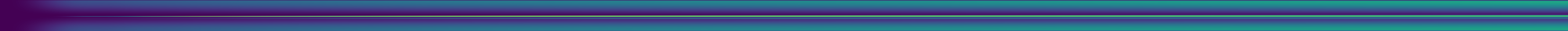};
    \end{axis}

    \begin{axis}[
        axis line style={draw=none},
        width=\textwidth/3,
        title={d: \(k = 10^{15};\ p = 2\)},
        name=plot3,
        at=(plot0.below south west), anchor=above north west,
        xlabel={Time \(t\) / \textmu s},
        ylabel={Position \(x\) / mm},
        xmin=0, xmax=499,
        ymin=0, ymax=1,
        title style={yshift=-7},
        tick align=outside,
        tick pos=left,
    ]
        \addplot graphics [xmin=0, xmax=499, ymin=0, ymax=1] {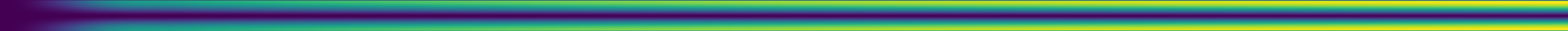};
    \end{axis}

    \begin{axis}[
        axis line style={draw=none},
        width=\textwidth/3,
        title={e: \(\al = 1.2;\ b = 10^{18}\)},
        name=plot4,
        at=(plot3.right of south east), anchor=left of south west,
        xlabel={Time \(t\) / \textmu s},
        xmin=0, xmax=499,
        ymin=0, ymax=1,
        yticklabels={,,},
        title style={yshift=-7},
        tick align=outside,
        tick pos=left,
    ]
        \addplot graphics [xmin=0, xmax=499, ymin=0, ymax=1] {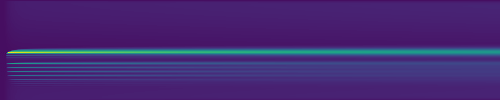};
    \end{axis}

    \begin{axis}[
        axis line style={draw=none},
        width=\textwidth/3,
        title={f: \(\al = 2;\ b = 10^{14}\)},
        name=plot5,
        at=(plot4.right of south east), anchor=left of south west,
        xlabel={Time \(t\) / \textmu s},
        xmin=0, xmax=499,
        ymin=0, ymax=1,
        yticklabels={,,},
        title style={yshift=-7},
        tick align=outside,
        tick pos=left,
        point meta min=0,
        point meta max=1,
        colorbar,
        colormap/viridis,
        colorbar style={
            axis line style={draw=none},
            ymin=0, ymax=1,
            ylabel={Norm. Temperature},
            tick align=outside,
            tick pos=right,
            at={(1.05,0.5)},
            anchor=west,
            width=0.04*\pgfkeysvalueof{/pgfplots/parent axis width},
        },
    ]
        \addplot graphics [xmin=0, xmax=499, ymin=0, ymax=1] {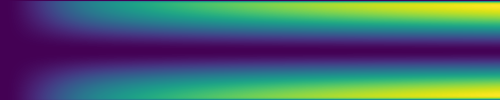};
    \end{axis}

\end{tikzpicture}
    \caption{Results for the thermal field for coupled thermo-acoustic simulations of a one-dimensional resonator with different values and different models for the dependence of the elasticity on temperature (a to d: power law, (\ref{equ:temp_dependence_pow}); e and f: convergent exponential expression, (\ref{equ:temp_dependence_exp})).}\label{fig:thermoacoustics_fields}
\end{figure}
For further analysis, the spatial temperature distributions for different values of \(k\) and \(p\) for the power law and \(\al\) and \(b\) convergent material behavior are examined.
The simulation results for the time-dependent temperature distribution for four different configurations are visualized in
Figure~\ref{fig:thermoacoustics_fields}.
It is immediately visible that qualitatively different results may arise, depending on the chosen parameter values.
Due to the mechanical boundary conditions, the temperature increase is expected to be maximal at the boundary of the system (\(x = 0\)\,mm and \(x = 1\)\,mm).
However, the thermal boundary conditions force the temperature to be zero at the same boundaries, resulting in the distributions shown in Figures~\ref{fig:thermoacoustics_fields}a, \ref{fig:thermoacoustics_fields}d or \ref{fig:thermoacoustics_fields}f, with local maxima close to but not at the spatial boundaries.
Aside from these results, qualitatively different temperature distributions (Figures~\ref{fig:thermoacoustics_fields}b and~\ref{fig:thermoacoustics_fields}e) are also observed.
In these cases, the temperature field forms a number of distinct, small hot spots along the spatial axis.
The number and periodicity of these hot spots depends on the choice of the parameters.
There are also observations, which show a superposition of the expected behavior and hot spots (Figure
\ref{fig:thermoacoustics_fields}c).

Hot spots are observed primarily for larger values of \(k\) when applying the power law, which is obvious when comparing
Figures~\ref{fig:thermoacoustics_fields}a and~\ref{fig:thermoacoustics_fields}b.
Increasing \(k\) further (e.g.\ to \(10^5\)) causes the simulation to destabilise.
Additionally, observations show that the temperature distribution will not show hot sports for values of \(p \leq 2\) for the power law.
For the convergent behavior, hot spots occur primarily for large values of \(b\), e.g.\ when the initial increase in elasticity \(C\) is large.
If the observed, qualitative differences in the spatial temperature distribution arise from effects of imminent numerical instability or if they also exist in physical systems is a subject to be explored in future research.

If it is found that the observed effects (formation of hot spots and unexpected rapid temperature increases) can occur in physical systems, even ideal cooling systems for acoustic resonators, such as high power piezoelectric actuators, may not be sufficient to keep such a system under stable operating conditions.
Because the cause for this behavior is an unfavourable temperature-dependence of the material, it may rule out certain material classes for an application in these systems.
Even if the observed artefacts are caused by numerical effect, further study of the coupled thermal and mechanical equation system is necessary to develop predictably stable simulation environments in the future.

\bigskip

{\bf Acknowledgement.} \quad
MW and LC acknowledge support of the {\em Deutsche Forschungsgemeinschaft} (No. 444955436).

\footnotesize 
  \setlength{\parskip}{0pt}
  \setlength{\itemsep}{0pt plus 0.2ex}


\begin{thebibliography}{10}

\bibitem{amann}
H.~Amann.
\newblock Nonhomogeneous linear and quasilinear elliptic and parabolic boundary
  value problems.
\newblock In {\em Function spaces, differential operators and nonlinear
  analysis ({F}riedrichroda, 1992)}, volume 133 of {\em Teubner-Texte Math.},
  pages 9--126. Teubner, Stuttgart, 1993.

\bibitem{bartels_roubicek}
S.~Bartels and T.~Roub\'{\i}\v{c}ek.
\newblock Thermoviscoplasticity at small strains.
\newblock {\em ZAMM Z. Angew. Math. Mech.}, 88(9):735--754, 2008.

\bibitem{bartels_roubicek_m2an11}
S.~Bartels and T.~Roub\'{\i}\v{c}ek.
\newblock Thermo-visco-elasticity with rate-independent plasticity in isotropic
  materials undergoing thermal expansion.
\newblock {\em ESAIM Math. Model. Numer. Anal.}, 45(3):477--504, 2011.

\bibitem{bartels_roubicek_numpde13}
S.~Bartels and T.~Roub\'{\i}\v{c}ek.
\newblock Numerical approaches to thermally coupled perfect plasticity.
\newblock {\em Numer. Methods Partial Differential Equations},
  29(6):1837--1863, 2013.

\bibitem{bies_cieslak}
P.~M. Bies and T.~Cieślak.
\newblock Time-asymptotics of a heated string.
\newblock 2024.
\newblock arXiv:2405.04310.

\bibitem{blanchard_guibe}
D.~Blanchard and O.~Guib\'{e}.
\newblock Existence of a solution for a nonlinear system in
  thermoviscoelasticity.
\newblock {\em Adv. Differential Equations}, 5(10-12):1221--1252, 2000.

\bibitem{Boley2012}
B.~Boley and J.~Weiner.
\newblock {\em Theory of Thermal Stresses}.
\newblock Dover Civil and Mechanical Engineering. Dover Publications, 2012.

\bibitem{bonetti_bonfanti}
E.~Bonetti and G.~Bonfanti.
\newblock Existence and uniqueness of the solution to a 3{D} thermoviscoelastic
  system.
\newblock {\em Electron. J. Differential Equations}, pages No. 50, 15, 2003.

\bibitem{chelminski_owczarek}
K.~Che\l{}mi\'{n}ski and S.~Owczarek.
\newblock Renormalised solutions in thermo-visco-plasticity for a
  {N}orton-{H}off type model. {P}art {II}: the limit case.
\newblock {\em Nonlinear Anal. Real World Appl.}, 31:643--660, 2016.

\bibitem{chelminski_owczarek_I}
K.~Che\l{}mi\'nski and S.~Owczarek.
\newblock Renormalized solutions in thermo-visco-plasticity for a
  {N}orton-{H}off type model. {P}art {I}: the truncated case.
\newblock {\em Nonlinear Anal. Real World Appl.}, 28:140--152, 2016.

\bibitem{cieslak_galic_muha}
T.~Cie\'{s}lak, M.~Gali\'{c}, and B.~Muha.
\newblock A model in one-dimensional thermoelasticity.
\newblock {\em Nonlinear Anal.}, 216:Paper No. 112703, 21, 2022.

\bibitem{cieslak_muha_trifunovic}
T.~Cie\'{s}lak, B.~Muha, and S.~a. Trifunovi\'{c}.
\newblock Global weak solutions in nonlinear 3{D} thermoelasticity.
\newblock {\em Calc. Var. Partial Differential Equations}, 63(1):Paper No. 26,
  36, 2024.

\bibitem{Claes2024}
L.~Claes and M.~Webersen.
\newblock pyfds 0.3.1 -- modular field simulation tool, 2024.

\bibitem{dafermos_hsiao_smooth}
C.~M. Dafermos and L.~Hsiao.
\newblock Global smooth thermomechanical processes in one-dimensional nonlinear
  thermoviscoelasticity.
\newblock {\em Nonlinear Anal.}, 6(5):435--454, 1982.

\bibitem{dafermos_hsiao}
C.~M. Dafermos and L.~Hsiao.
\newblock Development of singularities in solutions of the equations of
  nonlinear thermoelasticity.
\newblock {\em Quart. Appl. Math.}, 44(3):463--474, 1986.

\bibitem{Feldmann2021}
N.~Feldmann, V.~Schulze, L.~Claes, B.~Jurgelucks, L.~Meihost, A.~Walther, and
  B.~Henning.
\newblock Modelling damping in piezoceramics: A comparative study.
\newblock {\em tm - Technisches Messen}, 88(5):294--302, 2021.

\bibitem{friedman_book}
A.~Friedman.
\newblock {\em Partial differential equations}.
\newblock Holt, Rinehart and Winston, Inc., New York-Montreal, Que.-London,
  1969.

\bibitem{Friesen2023}
O.~Friesen, L.~Claes, C.~Scheidemann, N.~Feldmann, T.~Hemsel, and B.~Henning.
\newblock Estimation of temperature-dependent piezoelectric material parameters
  using ring-shaped specimens.
\newblock In {\em 2023 International Congress on Ultrasonics, Beijing, China},
  volume 2822, page 012125. IOP Publishing, 9 2024.

\bibitem{gawinecki03}
J.~Gawinecki.
\newblock Global existence of solutions for non-small data to non-linear
  spherically symmetric thermoviscoelasticity.
\newblock {\em Math. Methods Appl. Sci.}, 26(11):907--936, 2003.

\bibitem{gawinecki_zajaczkowski_cpaa}
J.~A. Gawinecki and W.~M. Zaj\c{a}czkowski.
\newblock Global regular solutions to two-dimensional thermoviscoelasticity.
\newblock {\em Commun. Pure Appl. Anal.}, 15(3):1009--1028, 2016.

\bibitem{gawinecki_zajaczkowski}
J.~A. Gawinecki and W.~M. Zaj\c{a}czkowski.
\newblock On regular solutions to two-dimensional thermoviscoelasticity.
\newblock {\em Appl. Math. (Warsaw)}, 43(2):207--233, 2016.

\bibitem{GutierrezLemini2014}
D.~Gutierrez-Lemini.
\newblock {\em {Engineering Viscoelasticity}}.
\newblock {Springer US}, Boston, MA and s.l., 2014.

\bibitem{henry}
D.~Henry.
\newblock {\em Geometric theory of semilinear parabolic equations}, volume 840
  of {\em Lecture Notes in Mathematics}.
\newblock Springer-Verlag, Berlin-New York, 1981.

\bibitem{higham}
N.~J. Higham.
\newblock {\em Functions of matrices}.
\newblock Society for Industrial and Applied Mathematics (SIAM), Philadelphia,
  PA, 2008.
\newblock Theory and computation.

\bibitem{horstmann_win}
D.~Horstmann and M.~Winkler.
\newblock Boundedness vs. blow-up in a chemotaxis system.
\newblock {\em J. Differential Equations}, 215(1):52--107, 2005.

\bibitem{jiang1990}
S.~Jiang.
\newblock Global existence of smooth solutions in one-dimensional nonlinear
  thermoelasticity.
\newblock {\em Proc. Roy. Soc. Edinburgh Sect. A}, 115(3-4):257--274, 1990.

\bibitem{jiang_racke90}
S.~Jiang and R.~Racke.
\newblock On some quasilinear hyperbolic-parabolic initial-boundary value
  problems.
\newblock {\em Math. Methods Appl. Sci.}, 12(4):315--339, 1990.

\bibitem{kim}
J.~U. Kim.
\newblock Global existence of solutions of the equations of one-dimensional
  thermoviscoelasticity with initial data in {$BV$} and {$L\sp{1}$}.
\newblock {\em Ann. Scuola Norm. Sup. Pisa Cl. Sci. (4)}, 10(3):357--427, 1983.

\bibitem{korn}
A.~Korn.
\newblock {\"U}ber einige {Ungleichungen}, welche in der {Theorie} der
  elastischen und elektrischen {Schwingungen} eine {Rolle} spielen.
\newblock Krak. {Anz}., 705-724 (1909)., 1909.

\bibitem{lankeit_win_NoDEA}
J.~Lankeit and M.~Winkler.
\newblock A generalized solution concept for the {K}eller-{S}egel system with
  logarithmic sensitivity: global solvability for large nonradial data.
\newblock {\em NoDEA Nonlinear Differential Equations Appl.}, 24(4):Paper No.
  49, 33, 2017.

\bibitem{Lesieutre1996}
G.~A. Lesieutre, L.~Fang, G.~H. Koopmann, S.~P. Pai, and S.~Yoshikawa.
\newblock Heat generation of a piezoceramic induced-strain actuator embedded in
  a glass/epoxy composite panel.
\newblock In I.~Chopra, editor, {\em Smart Structures and Materials 1996: Smart
  Structures and Integrated Systems}. SPIE, 5 1996.

\bibitem{mcintire}
R.~S. McIntire.
\newblock {\em A new technique for discussing the development of singularities
  in quasilinear hyperbolic PDE's with applications to a model problem in
  nonlinear theormoelasticity (blow-up, catastrophe, breakdown)}.
\newblock Brown University, 1985.

\bibitem{Meyers2008}
M.~A. Meyers and K.~K. Chawla.
\newblock {\em Mechanical behavior of materials}.
\newblock Cambridge University Press, Cambridge, England, 2 edition, 11 2008.

\bibitem{mielke_roubicek}
A.~Mielke and T.~Roub\'{\i}\v{c}ek.
\newblock Thermoviscoelasticity in {K}elvin-{V}oigt rheology at large strains.
\newblock {\em Arch. Ration. Mech. Anal.}, 238(1):1--45, 2020.

\bibitem{neff_pauly_witsch}
P.~Neff, D.~Pauly, and K.-J. Witsch.
\newblock Poincar{\'e} meets {Korn} via {Maxwell}: extending {Korn}'s first
  inequality to incompatible tensor fields.
\newblock {\em J. Differ. Equations}, 258(4):1267--1302, 2015.

\bibitem{Nye1985}
J.~Nye.
\newblock {\em Physical Properties of Crystals: Their Representation by Tensors
  and Matrices}.
\newblock Oxford science publications. Clarendon Press, 1985.

\bibitem{owczarek_wielgos}
S.~Owczarek and K.~Wielgos.
\newblock On a thermo-visco-elastic model with nonlinear damping forces and
  {$L^1$} temperature data.
\newblock {\em Math. Methods Appl. Sci.}, 46(9):9966--9999, 2023.

\bibitem{paoli_petrov}
L.~Paoli and A.~Petrov.
\newblock Global existence result for thermoviscoelastic problems with
  hysteresis.
\newblock {\em Nonlinear Anal. Real World Appl.}, 13(2):524--542, 2012.

\bibitem{paoli_petrov_gamm}
L.~Paoli and A.~Petrov.
\newblock Thermodynamics of multiphase problems in viscoelasticity.
\newblock {\em GAMM-Mitt.}, 35(1):75--90, 2012.

\bibitem{pawlow_zajaczkowski_cpaa17}
I.~Paw\l{}ow and W.~M. Zaj\c{a}czkowski.
\newblock Global regular solutions to three-dimensional thermo-visco-elasticity
  with nonlinear temperature-dependent specific heat.
\newblock {\em Commun. Pure Appl. Anal.}, 16(4):1331--1371, 2017.

\bibitem{PICeramic2023}
PI Ceramic GmbH.
\newblock {\em Material Data -- Specific parameters of the standard materials},
  2023.

\bibitem{qin}
Y.~Qin.
\newblock Global existence and asymptotic behaviour of the solution to the
  system in one-dimensional nonlinear thermoviscoelasticity.
\newblock {\em Quart. Appl. Math.}, 59(1):113--142, 2001.

\bibitem{racke}
R.~Racke.
\newblock Initial boundary value problems in one-dimensional nonlinear
  thermoelasticity.
\newblock {\em Math. Methods Appl. Sci.}, 10(5):517--529, 1988.

\bibitem{racke_bu}
R.~Racke.
\newblock Blow-up in nonlinear three-dimensional thermoelasticity.
\newblock {\em Math. Methods Appl. Sci.}, 12(3):267--273, 1990.

\bibitem{racke90}
R.~Racke.
\newblock On the {C}auchy problem in nonlinear {$3$}-d thermoelasticity.
\newblock {\em Math. Z.}, 203(4):649--682, 1990.

\bibitem{racke_zheng}
R.~Racke and S.~Zheng.
\newblock Global existence and asymptotic behavior in nonlinear
  thermoviscoelasticity.
\newblock {\em J. Differential Equations}, 134(1):46--67, 1997.

\bibitem{Remis2000}
R.~F. Remis.
\newblock On the stability of the finite-difference time-domain method.
\newblock {\em Journal of Computational Physics}, 163(1):249--261, 9 2000.

\bibitem{rossi_roubicek}
R.~Rossi and T.~Roub\'{\i}\v{c}ek.
\newblock Thermodynamics and analysis of rate-independent adhesive contact at
  small strains.
\newblock {\em Nonlinear Anal.}, 74(10):3159--3190, 2011.

\bibitem{rossi_roubicek_interfaces13}
R.~Rossi and T.~Roub\'{\i}\v{c}ek.
\newblock Adhesive contact delaminating at mixed mode, its thermodynamics and
  analysis.
\newblock {\em Interfaces Free Bound.}, 15(1):1--37, 2013.

\bibitem{roubicek}
T.~Roub\'{\i}\v{c}ek.
\newblock Thermo-visco-elasticity at small strains with {$L^1$}-data.
\newblock {\em Quart. Appl. Math.}, 67(1):47--71, 2009.

\bibitem{roubicek_SIMA10}
T.~Roub\'{\i}\v{c}ek.
\newblock Thermodynamics of rate-independent processes in viscous solids at
  small strains.
\newblock {\em SIAM J. Math. Anal.}, 42(1):256--297, 2010.

\bibitem{roubicek_nodea13}
T.~Roub\'{\i}\v{c}ek.
\newblock Nonlinearly coupled thermo-visco-elasticity.
\newblock {\em NoDEA Nonlinear Differential Equations Appl.}, 20(3):1243--1275,
  2013.

\bibitem{roubicek_dcdss13}
T.~Roub\'{\i}\v{c}ek.
\newblock Thermodynamics of perfect plasticity.
\newblock {\em Discrete Contin. Dyn. Syst. Ser. S}, 6(1):193--214, 2013.

\bibitem{Rudenko1977}
O.~V. Rudenko, S.~I. Solujan, and R.~T. Beyer.
\newblock {\em Theoretical foundations of nonlinear acoustics}.
\newblock Studies in Soviet science. {Consultants Bureau}, New York, 1977.

\bibitem{Rupitsch2019}
S.~J. Rupitsch.
\newblock {\em Piezoelectric sensors and actuators}.
\newblock Springer, 2019.

\bibitem{shibata}
Y.~Shibata.
\newblock Global in time existence of small solutions of nonlinear
  thermoviscoelastic equations.
\newblock {\em Math. Methods Appl. Sci.}, 18(11):871--895, 1995.

\bibitem{simon}
J.~Simon.
\newblock Compact sets in the space {$L^p(0,T;B)$}.
\newblock {\em Ann. Mat. Pura Appl. (4)}, 146:65--96, 1987.

\bibitem{slemrod}
M.~Slemrod.
\newblock Global existence, uniqueness, and asymptotic stability of classical
  smooth solutions in one-dimensional nonlinear thermoelasticity.
\newblock {\em Arch. Rational Mech. Anal.}, 76(2):97--133, 1981.

\bibitem{Tauchert1967}
T.~R. Tauchert.
\newblock {Heat generation in a viscoelastic solid}.
\newblock {\em {Acta Mechanica}}, 3(4):385--396, 1967.

\bibitem{truesdell_cauchy}
C.~A. Truesdell.
\newblock Cauchy and the modern mechanics of continua.
\newblock volume~45, pages 5--24. 1992.
\newblock \'{E}tudes sur Cauchy (1789--1857).

\bibitem{Wellendorf2023}
A.~Wellendorf, L.~von Damnitz, A.~Nuri, D.~Anders, and S.~Trampnau.
\newblock Determination of the temperature-dependent resonance behavior of
  ultrasonic transducers using the finite-element method.
\newblock {\em Journal of Vibration Engineering \& Technologies}, pages 1--14,
  2 2023.

\bibitem{win_M3AS_scrounge}
M.~Winkler.
\newblock Global generalized solutions to a multi-dimensional doubly tactic
  resource consumption model accounting for social interactions.
\newblock {\em Math. Models Methods Appl. Sci.}, 29(3):373--418, 2019.

\bibitem{Yee1966}
K.~Yee.
\newblock Numerical solution of initial boundary value problems involving
  {M}axwell{\textquotesingle}s equations in isotropic media.
\newblock {\em {IEEE} Transactions on Antennas and Propagation},
  14(3):302--307, 5 1966.

\bibitem{yoshikawa_pawlow_zajaczkowski_SIMA}
S.~Yoshikawa, I.~Paw\l{}ow, and W.~M. Zaj\c{a}czkowski.
\newblock Quasi-linear thermoelasticity system arising in shape memory
  materials.
\newblock {\em SIAM J. Math. Anal.}, 38(6):1733--1759, 2007.

\bibitem{yoshikawa_pawlow_zajaczkowski}
S.~Yoshikawa, I.~Paw\l{}ow, and W.~M. Zaj\c{a}czkowski.
\newblock A quasilinear thermoviscoelastic system for shape memory alloys with
  temperature dependent specific heat.
\newblock {\em Commun. Pure Appl. Anal.}, 8(3):1093--1115, 2009.

\bibitem{Zheng1996}
J.~Zheng, S.~Takahashi, S.~Yoshikawa, K.~Uchino, and J.~W.~C. de~Vries.
\newblock Heat generation in multilayer piezoelectric actuators.
\newblock {\em Journal of the American Ceramic Society}, 79(12):3193--3198,
  1996.

\end{thebibliography}
\end{document}